\documentstyle[amstex]{amsart}
\newtheorem{Theorem}{Theorem}[section]
\newtheorem{Lemma}{Lemma}[section]
\newtheorem{Remark}{Remark}[section]

\newtheorem{Proposition}{Proposition}[section]
\newtheorem{Example}{Example}[section]
\newtheorem{Definition}{Definition}[section]

\newcommand\Pic{\mathop{\rm Pic}\nolimits}
\newcommand\Spec{\mathop{\rm Spec}\nolimits}

\makeatletter
\begin{document}

\title{Classification of Okamoto--Painlev\'e Pairs}

\author{Masa-Hiko Saito} 
\thanks{Partly supported by Grant-in Aid
for Scientific Research (B-09440015), the Ministry of
Education, Science and Culture, Japan.}
\author{Taro Takebe}
\address{Department of Mathematics, Faculty of Science, 
Kobe Unversity, Kobe, Rokko, 657-8501, Japan.}
\email{mhsaito@math.kobe-u.ac.jp}
\email{takebe@math.kobe-u.ac.jp}

\begin{abstract}
In this paper, we introduce the notion of an Okamoto--Painlev\'e pair 
$(S, Y)$ which consists of a compact smooth complex surface $S$ and an 
effective divisor $Y$ on $S$ satisfying certain conditions. 
Though spaces of initial values of Painlev\'e equations 
introduced by K. Okamoto give examples of Okamoto--Painleve pairs, 
we find a new example of Okamoto--Painlev\'e pairs not listed in \cite{Oka}.
We will give the complete classification of Okamoto--Painlev\'e pairs.
\end{abstract}  
\maketitle
\addtocounter{section}{-1}

\section{Introduction}
In this paper, we will introduce the notion of an Okamoto--Painlev\'e pair
$(S,Y)$, which is defined as follows:
\begin{Definition}
\label{def:op}
{\rm (Cf. Definition \ref{def:op1})}.
Let $S$ be a compact smooth complex surface and $Y = \sum_{i=1}^r a_i Y_i$ 
an effective divisor on $S$.  
We say that
a pair $(S, Y)$ is an Okamoto--Painlev\'e pair if it satisfies the 
following conditions:

{\rm (i)} There exists a meromorphic $2$-form $\omega$ on $S$ 
such that $(\omega) = -Y$, that is, $\omega$ has the  pole 
divisor $Y$ {\rm (}counting multiplicities{\rm )} and has no zero outside 
$Y$. 

{\rm (ii)} For all $i\;(1 \leq i \leq r)$, 
$Y \cdot Y_i = \deg [Y]|_{Y_i} = 0.$ 

{\rm (iii)} Let us set  $D :=  Y_{red} = \sum_{i=1}^{r} Y_i$.
Then $S - D$ contains ${\bf C}^2$ as a Zariski open set.

{\rm (iv)} Set $F = S - {\bf C}^2$ where ${\bf C}^2$ is the same Zariski 
open set  as in {\rm (ii)}. 
Then $F$ is a {\rm (}reduced{\rm )} divisor with normal crossings.
\end{Definition}
Historically, Okamoto \cite{Oka} introduced 
{\em the space $M_{J}(t)$ of initial values}
for each Painlev\'e equation of type $P_J\;(J=I,\ldots ,VI)$ 
with the time parameter $t$, which is a noncompact complex surface.
We can obtain a nice compactification $\overline{M}_{J}(t)$ of $M_{J}(t)$
so that $(\overline{M}_{J}(t),D_{J}(t))$ becomes an Okamoto--Painlev\'e pair,
where $D_{J}(t)=\overline{M}_{J}(t) -M_{J}(t)$
(cf.\ \cite{ST} and \cite{MMT}).
Conversely, for each Okamoto--Painlev\'e pair one can 
associate a Hamiltonian system via the deformation theory of pairs,
and such a Hamiltonian system is equivalent to a differential equation
of Painlev\'e type
(cf.\ \cite{Sa--T}).

The main purpose of this paper is the classification of Okamoto--Painlev\'e 
pairs.
We will classify Okamoto--Painlev\'e pairs $(S,Y)$ into seven
types according to the configuration of the divisor $Y\in |-K_S|$.
(See Theorem \ref{clop}).
One can easily see that the configuration of the divisor $Y$ is same as one of
singular fibers of elliptic surfaces \cite{Kod1}. 
According to the Kodaira's classification, the seven configurations 
can be denoted by
$II^{\ast}$, $III^{\ast}$, $IV^{\ast}$, $I_{3}^{\ast}$, $I_{2}^{\ast}$,
$I_{1}^{\ast}$, $I_{0}^{\ast}$,
or in the notation of the dual graph of $Y$ 
they can be denoted by
\[\widetilde{E}_8 , \widetilde{E}_7 , \widetilde{E}_6 , \widetilde{D}_7 , 
\widetilde{D}_6 , \widetilde{D}_5 , \widetilde{D}_4,\]
respectively.

It is known that
there is a one-to-one correspondence between Okamoto--Painlev\'e pairs
of the six types
$\widetilde{E}_8 , \widetilde{E}_7 , \widetilde{E}_6 , \widetilde{D}_7 , 
\widetilde{D}_6 , \widetilde{D}_5 , \widetilde{D}_4$
and Painlev\'e equations of six types
$P_{I}$, $P_{II}$, $P_{III}$, $P_{IV}$, $P_{V}$, $P_{VI}$, respectively
(cf.\ \cite{Oka}, \cite{MMT} and \cite{ST}).
Though an Okamoto--Painlev\'e pair of type $\widetilde{D}_7$ did not appear in 
the list of Okamoto \cite{Oka},
we can obtain a Hamiltonian system associated to the pairs.
Recently, Okamoto informed us that he found Painlev\'e equations for an
Okamoto--Painlev\'e pair of type $\widetilde{D}_7$ 
which are special cases of $P_{III}$.
Therefore, we denote this Painlev\'e equation by $P_{III}^{\ast}$.
In \S\ref{cstrD7}, we will construct an Okamoto--Painlev\'e pair of type
$\widetilde{D}_7$ by blowing-up of 
\({\bf F}_0={\bf P}^{1}\times{\bf P}^{1}\)
with centers on the anti-canonical divisor $-K_{{\bf F}_{0}}$. 

Let us discuss main ideas of classfying Okamoto--Painlev\'e pairs.
For an Okamoto--Painlev\'e pair $(S,Y)$, the surface $S$ is the 
compactification of ${\bf C}^{2}$, namely $S$ includes ${\bf C}^{2}$
as a Zariski open set. Let \(F=S-{\bf C}^{2}\).
By the definition of Okamoto--Painlev\'e pairs, $F$ is a 
normal crossing divisor. 
Each component of $Y_{\it red}$ is also a component of $F$, namely,
\[Y_{\it red} = \sum_{i=1}^{r} Y_i \subset \sum_{i=1}^{l} F_i =F,\]
hence $Y_{red}$ is also a normal crossing divisor.
One can easily see that the divisor $Y$ has the same configuration of
a singular fiber of an elliptic surface, 
and the number of irreducible components of $F$ is 10.
Therefore, the number of the irreducible components of $Y_{red}$ is less than 
or equal to 10.
Moreover, the configuration of $F$ must be a tree, so is $Y_{red}$.
These arguments show that the configuration of $Y$ is one of the types
$\widetilde{E}_{r-1}$ for $r=9,8,7$ and $\widetilde{D}_{r-1}$ for $r=5,\ldots,
10$.

Now we will consider the classification of the configurations
of normal clossing divisors $F$.
The divisor $F$ can be obtained by adding some components to $Y_{\it red}$.
We call an irreducible component of $F-Y_{red}$ an  {\em additional component}.
We will classify all of configurations of $F$, 
and this gives the complete classification of Okamoto--Painlev\'e pairs
and also the configurations of $F$.

The complete list of all normal crossing divisors $F$ is given in \S\ref{s31}.
Each of those is a divisor with ten components including one of the seven
Kodaira types.
All configurations of $F$ in the list can be transformed into the
anti-canonical divisors of \({\bf P}^{2},{\bf F}_{0}\) and \({\bf F}_{2}\) 
by blowings-up and blowings-down.

Let us discuss the relation between our results and results of Sakai in 
\cite{Sak}.
Sakai defined the notion of {\em generalized Halphen surfaces}.
By definition, a generalized Halphen surface $S$ is a compact complex surface
satisfying the condition of (i) and (ii) of Definition \ref{def:op}.
He related generalized Halphen surfaces to the discrete Painlev\'e equations
via Cremona transformations.
The Painlev\'e differential equations can be obtained as limits of discrete
Painlev\'e equations.
Most essential extra conditions of an Okamoto--Painlev\'e pair $(S,Y)$ 
are that $S$ is a compactification of ${\bf C}^2$ and $F$ is a normal crossing
divisor.
Our classification gives more direct correspondences between 
Okamoto--Painlev\'e pairs and Painlev\'e differential equations (cf.\ 
\cite{Sa--T}).

\section{Preliminary}

Let $S$ be a compact complex surface and let $-K_{S}$ denote the 
anti-canonical divisor class of $S$.
For any divisor $D$ on $S$, we denote by $[D]$ linear equivarence class of $D$,
and by $|D|=|[D]|$ the linear space of effective divisors $C$, such that 
$C\sim D$. 
Here ``$\sim$'' means the linear equivalence of divisors.
In this paper, we often identify the divisor class $[D]$ with 
the isomorphism class of the corresponding line bundle $[D]$ or 
the correponding invertible sheaf ${\cal O}_S(D)$.

Let us assume that there exists a normal crossing effective divisor
$Y\in |-K_S|$.
Moreover, assume that every irreducible component of $Y$ is a smooth rational
curve.
Set $ Y =\sum_{i=1}^{n} m_{i}C_{i}$.
Consider the blowing-up \({\pi}:\widetilde{S}\rightarrow S\) with
center of $P\in C_{i_{0}}$, 
and let $C_{i}'$ be the strict transform of $C_{i}$ by $\pi$, and $E$
the exceptional curve of $\pi$.
Note that $m_{i}\geq 0$.
By a standard calculation, we can obtain the following 
lemma.

\begin{Lemma}
If $p\in C_{i_{0}}\setminus \bigcup_{i\neq i_{0}} (C_{i_{0}}\cap C_{i})$,
then 
\[-K_{\widetilde{S}}=[\sum_{i=1}^{n} m_{i}C_{i}' +(m_{i_{0}}-1)E],\]
\[(C_{i}')^{2}=(C_{i})^{2} \quad (i\neq i_{0}),\]
and
\[(C_{i_{0}}')^{2}=(C_{i_{0}})^{2}-1.\]
If $p\in C_{i_{0}}\cap C_{i_{1}}$ for some $i_{1}\;(\neq i_{0})$,
then 
\[-K_{\widetilde{S}}=[\sum_{i=1}^{n} m_{i}C_{i}' 
+(m_{i_{0}}+m_{i_{1}}-1)E],\]
\[(C_{i}')^{2}=(C_{i})^{2} \quad (i\neq i_{0},i_{1}),\]
and
\[(C_{i}')^{2}=(C_{i})^{2}-1 \quad (i=i_{0},i_{1}).\vspace{0.5cm}\]
\end{Lemma}

Let ${\bf F}_{n}$ denote the Hirzebruch surface
${\bf P}({\cal O}_{{\bf P}^{1}}\oplus{\cal O}_{{\bf P}^{1}}(n)).$
We denote by $f$ the class of fiber of ${\bf F}_{n}$, $s_{0}$ the class of 
the minimal section of ${\bf F}_{n}$ respectively.
We also set $s_{\infty}=s_{0}+nf$. 
Then, the self-intersection numbers are given by 
$$
f^{2}=0,\quad (s_{0})^{2}=-n,\quad (s_{\infty})^{2}=n.
$$
It is well-known that ${\bf F}_{n}$ can be obtained by performing 
several blowings-up and blowing-downs to ${\bf P}^{2}$.

Now we consider the anti-canonical divisor class of ${\bf P}^{2}$ and
${\bf F}_{n}$.
Let $h$ denote the class of lines in ${\bf P}^{2}$.
Then, they are given as
\begin{eqnarray*}
-K_{{\bf P}^{2}} &=& 3h,\\
-K_{{\bf F}_{n}} &=& 2s_{0}+(n+2)f= 2s_{\infty}-(n-2)f.
\end{eqnarray*}
In this paper, we often use the formula $-K_{{\bf F}_{0}}=2s_{0}+2f$
and $-K_{{\bf F}_{2}}=2s_{\infty}$.
Note that 
$-K_{{\bf F}_{0}}$ and $-K_{{\bf F}_{2}}$
have complements in ${\bf F}_{0}$ and ${\bf F}_{2}$ respectively 
which contain ${\bf C}^{2}$ as Zariski open sets.

\section{Okamoto--Painlev\'e Pairs}
Now let us give the definition of Okamoto--Painlev\'e pairs.

\begin{Definition}
\label{def:op1}
Let $S$ be a compact smooth complex surface and $Y = \sum_{i=1}^r a_i Y_i$ 
an effective divisor on $S$.  We say that 
a pair $(S, Y)$ is an Okamoto--Painlev\'e pair if it satisfies the 
following conditions:

{\rm (i)} There exists a meromorphic $2$-form $\omega$ on $S$ 
such that $(\omega) = -Y$, that is, $\omega$ has the  pole 
divisor $Y$ {\rm (}counting multiplicities{\rm )} and has no zero outside 
$Y$. 

{\rm (ii)} For all $i\;(1 \leq i \leq r)$, 
$Y \cdot Y_i = \deg [Y]|_{Y_i} = 0.$ 

{\rm (iii)} Let us set  $D :=  Y_{red} = \sum_{i=1}^{r} Y_i$.
Then $S - D$ contains ${\bf C}^2$ as a Zariski open set.

{\rm (iv)} Set $F = S - {\bf C}^2$ where ${\bf C}^2$ is the same Zariski 
open set  as in {\rm (ii)}. 
Then $F$ is a {\rm (}reduced{\rm )} divisor with normal crossings.
\end{Definition}

\begin{Remark}
{\rm (i)} The meromorphic $2$-form $\omega$ as above
induces a holomorphic symplectic structure on $S-D$.

{\rm (ii)} Let $K_S$ denote the canonical divisor class of $S$.  
The condition {\rm (i)} means that $K_S = [-Y]$ or $-K_S=[Y]$.

{\rm (iii)} The condition {\rm (v)} implies that the reduced part 
$D= Y_{red}$ of $Y$ 
is also a divisor with normal crossings.
 
{\rm (iv)} The spaces of initial values of 
Painlev\'e equations can be written as $S - D$ 
for some Okamoto--Painlev\'e pair $(S, Y)$.
{\rm (See \cite{Oka}, \cite{ST} and \cite{MMT})}.

\end{Remark}

The following is the main theorem of this paper.
\begin{Theorem}\label{clop}
Let $(S, Y)$ be an Okamoto--Painlev\'e pair. 
Then we have the following assertions.

{\rm (i)} The surface $S$ is a projective rational surface.

{\rm (ii)} The configuration of $Y$ counting
multiplicities  is in the list of 
Kodaira's classification of singular fibers of 
elliptic surfaces.  More precisely, they are 
given by one of the following lists.

{\rm(iii)} All Okamoto--Painlev\'e pairs can be obtained by blowings-up
and blowings-down of $({\bf P}^2, 3H)$, where $H$ is a line in 
${\bf P}^2$.

\vspace{0.2cm}

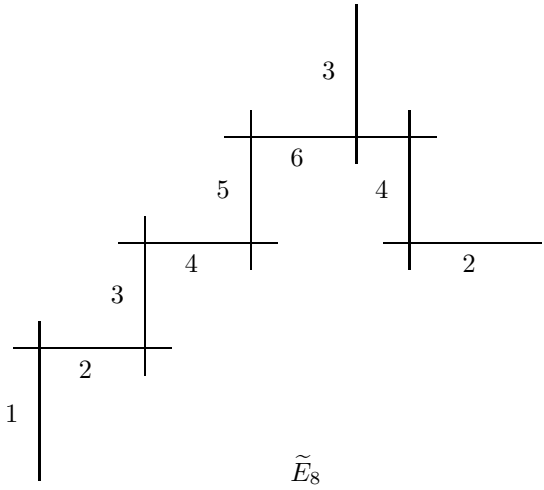
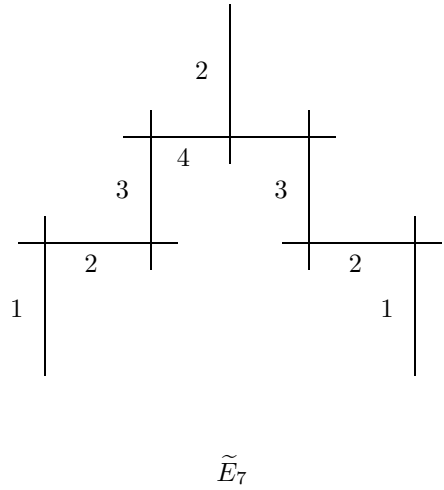
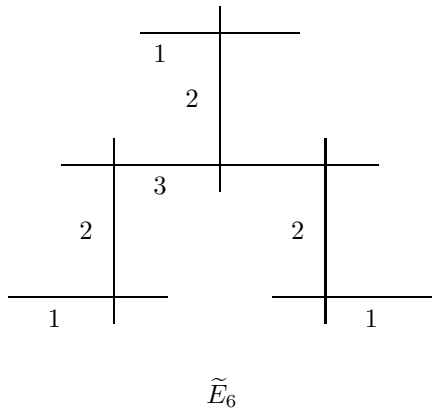
\begin{center}

\begin{tabular}{||c||c|c|c|c|c|c|c|c||} \hline
    &   & & & & & &  \\
$Y$ & $\widetilde{E}_8$ & $\widetilde{E}_7$ & $\widetilde{D}_7$ & 
$\widetilde{D}_6$ & $\widetilde{E}_6$ &$\widetilde{D}_5$ & $\widetilde{D}_4$ \\
    &   & & & & & &  \\ \hline
        &   & & & & & &  \\
Kodaira's notation & $II^*$ & $III^*$ & $I_{3}^*$ & $I_{2}^*$ & 
$IV^*$ & $I_{1}^*$ & $I_{0}^*$ \\ 
    &   & & & & & &  \\ \hline
        &   & & & & & &  \\ 
Painlev\'{e} equation & $P_{I}$  & $P_{II}$ & $P_{III}^{\ast}$ & $P_{III}$
& $P_{IV}$ 
& $P_{V}$  &  $P_{VI}$ \\
    &   & & & & & &  \\ \hline
\end{tabular}

\end{center}
\vspace{1.2cm}

\begin{picture}(400,200)
\put(-10,0){\line(0,1){60}}
\put(-20,50){\line(1,0){60}}
\put(30,40){\line(0,1){60}}
\put(20,90){\line(1,0){60}}
\put(70,80){\line(0,1){60}}
\put(60,130){\line(1,0){80}}
\put(110,120){\line(0,1){60}}
\put(130,80){\line(0,1){60}}
\put(120,90){\line(1,0){60}}

\put(-23,22){$1$}
\put(5,39){$2$}
\put(17,67){$3$}
\put(45,79){$4$}
\put(57,107){$5$}
\put(85,119){$6$}
\put(97,152){$3$}
\put(117,107){$4$}
\put(150,79){$2$}

\put(85,0){$\widetilde{E}_{8}$}
\put(250,40){\line(0,1){60}}
\put(240,90){\line(1,0){60}}
\put(290,80){\line(0,1){60}}
\put(280,130){\line(1,0){80}}
\put(320,120){\line(0,1){60}}
\put(350,80){\line(0,1){60}}
\put(340,90){\line(1,0){60}}
\put(390,40){\line(0,1){60}}

\put(237,62){$1$}
\put(265,79){$2$}
\put(277,107){$3$}
\put(300,119){$4$}
\put(307,152){$2$}
\put(337,107){$3$}
\put(365,79){$2$}
\put(377,62){$1$}

\put(315,0){$\widetilde{E}_{7}$}
\end{picture}

\vspace{1cm}

\begin{picture}(400,150)
\put(30,80){\line(1,0){120}}
\put(90,70){\line(0,1){70}}
\put(60,130){\line(1,0){60}}
\put(50,20){\line(0,1){70}}
\put(10,30){\line(1,0){60}}
\put(130,20){\line(0,1){70}}
\put(110,30){\line(1,0){60}}

\put(65,69){$3$}
\put(77,102){$2$}
\put(65,119){$1$}
\put(37,52){$2$}
\put(25,19){$1$}
\put(117,52){$2$}
\put(145,19){$1$}

\put(85,-10){$\widetilde{E}_{6}$}
\end{picture}

\vspace{1cm}

\begin{picture}(400,200)
\put(-10,70){\line(1,0){100}}
\put(10,20){\line(0,1){60}}
\put(70,20){\line(0,1){60}}
\put(50,60){\line(0,1){90}}
\put(30,130){\line(1,0){90}}
\put(110,90){\line(0,1){100}}
\put(100,110){\line(1,0){60}}
\put(100,170){\line(1,0){60}}

\put(25,59){$2$}
\put(-6,42){$1$}
\put(54,42){$1$}
\put(34,97){$2$}
\put(75,119){$2$}
\put(94,147){$2$}
\put(130,99){$1$}
\put(130,159){$1$}
\put(100,0){$\widetilde{D}_{7}$}
\put(260,100){\line(1,0){110}}
\put(280,40){\line(0,1){120}}
\put(230,60){\line(1,0){60}}
\put(230,140){\line(1,0){60}}
\put(350,40){\line(0,1){120}}
\put(340,60){\line(1,0){60}}
\put(340,140){\line(1,0){60}}

\put(310,89){$2$}
\put(264,77){$2$}
\put(250,49){$1$}
\put(250,129){$1$}
\put(334,77){$2$}
\put(370,49){$1$}
\put(370,129){$1$}
\put(315,0){$\widetilde{D}_{6}$}
\end{picture}

\vspace{1cm}

\begin{picture}(400,150)
\put(0,50){\line(1,0){120}}
\put(20,0){\line(0,1){60}}
\put(70,0){\line(0,1){60}}
\put(100,30){\line(0,1){120}}
\put(90,80){\line(1,0){60}}
\put(90,130){\line(1,0){60}}

\put(40,39){$2$}
\put(4,22){$1$}
\put(54,22){$1$}
\put(84,102){$2$}
\put(120,69){$1$}
\put(120,119){$1$}
\put(100,-20){$\widetilde{D}_{5}$}

\put(240,70){\line(1,0){160}}
\put(260,30){\line(0,1){80}}
\put(300,30){\line(0,1){80}}
\put(340,30){\line(0,1){80}}
\put(380,30){\line(0,1){80}}

\put(315,59){$2$}
\put(247,87){$1$}
\put(287,87){$1$}
\put(327,87){$1$}
\put(367,87){$1$}
\put(315,-20){$\widetilde{D}_{4}$}
\put(170,-60){{\bf Figure 1.}}

\end{picture}
\end{Theorem}
\vspace{3cm}
In Figure 1, each line denotes ${\bf P}^{1}$ whose self-intersection number
is equal to $-2$, 
and the number next to each line denotes the multiplicity of corresponding component in $-K_{S}$.\\

Let us start our proof with the following easy lemmas.
\begin{Lemma}
\label{lem1}
Let $(S, Y)$ be an Okamoto--Painlev\'e pair, then we have:\\
{\rm (i)} $Y_{i}\cong {\bf P}^{1}$ for $i\; (1\leq i \leq r)$.\\
{\rm (ii)} The following conditions are equivarent for $i\;
(1\leq i \leq r)$.

{\rm (a)} $Y\cdot Y_{i} =0$.

{\rm (b)} $(Y_{i})^{2} =-2$.
\end{Lemma}
\noindent
{\em Proof.} 
{\rm (i)} It is sufficient to show that the component $F_{i}$ 
of $F$ is isomorphic to ${\bf P}^{1}$. 
We will show that $H^{1}(F_{i},{\bf Z})=0$.\\
We have the following isomorphisms of cohomology groups with
{\bf Z}-coefficients by Poincar\'e duality.
\[ \begin{array}{ccl}
H^{i}(S,F;{\bf Z}) & {\cong} & H^{i}_{c}({\bf C}^{2},{\bf Z})\\
                   & {\cong} & H_{4-i}({\bf C}^{2},{\bf Z})\\
                   & {\cong} & \left\{ \begin{array}{ll}
                                       0 & i=0,1,2,3\\
                                       {\bf Z} & i=4
                                       \end{array}
                               \right.

\end{array}\]
On the other hand, consider the long exact sequence of cohomology groups
for pair $(S,F)$
\[\begin{array}{cclclcl}
0 & \rightarrow & H^{0}(F,{\bf Z}) & \rightarrow & H^{0}(S,{\bf Z}) &
\rightarrow & H^{0}(S,F;{\bf Z})=0\\
& \rightarrow & H^{1}(F,{\bf Z}) & \rightarrow & H^{1}(S,{\bf Z}) &
\rightarrow & H^{1}(S,F;{\bf Z})=0\\
& \rightarrow & H^{2}(F,{\bf Z}) & \rightarrow & H^{2}(S,{\bf Z}) &
\rightarrow & H^{2}(S,F;{\bf Z})=0\\
& \rightarrow & H^{3}(F,{\bf Z}) & \rightarrow & H^{3}(S,{\bf Z}) &
\rightarrow & H^{3}(S,F;{\bf Z})=0\\
& \rightarrow & H^{4}(F,{\bf Z}) & \rightarrow & H^{4}(S,{\bf Z}) &
\rightarrow & H^{4}(S,F;{\bf Z})\cong {\bf Z}\\
& \rightarrow & \ldots. & & & &
\end{array}\]
Then we have 
\begin{equation}\label{a1}
H^{i}(S,{\bf Z}) \cong H^{i}(F,{\bf Z})\quad (i=0,1,2,3).
\end{equation}
Especially,
\[H^{3}(S,{\bf Z}) \cong H^{3}(F,{\bf Z}) =0.\]
By Poincar\'e duality again, we have
\begin{equation}\label{a2}
H^{1}(S,{\bf Z})\cong H^{3}(S,{\bf Z})=0.
\end{equation}
Now we consider an irreducible component $F_{i}$ of $F$.
Let \(F'=\bigcup_{k\neq i} F_{k}\).
Then $F_{i}\cap F'$ consists of a finite set of points.
And we have the Mayer--Vietoris exact sequence
\begin{equation}\label{a3}
     \begin{array}{cclclcl}
      0 & \rightarrow & H^{0}(F,{\bf Z}) & \rightarrow &
      H^{0}(F_{i},{\bf Z})\oplus H^{0}(F',{\bf Z}) & \rightarrow &
      H^{0}(F_{i}\cap F',{\bf Z})\\
      & \rightarrow & H^{1}(F,{\bf Z}) & \rightarrow &
      H^{1}(F_{i},{\bf Z})\oplus H^{1}(F',{\bf Z}) & \rightarrow &
      H^{1}(F_{i}\cap F',{\bf Z})\\
      & \rightarrow & \ldots. & & & &
     \end{array}
\end{equation}
By (\ref{a1}) and (\ref{a2}), we have $H^{1}(F,{\bf Z})= 0$.
Since \(H^{1}(F_{i}\cap F',{\bf Z})=0\), we have 
\(H^{1}(F_{i},{\bf Z})\oplus H^{1}(F',{\bf Z})=0\).
Therefore, $H^{1}(F_{i},{\bf Z})=0$.\\
If the irreducible component $F_{i}$ is a singular nodal curve, 
then $H^{1}(F_{i},{\bf Z})\neq 0$.
So $F_{i}$ is nonsingular, namely $F_{i}\cong {\bf P}^{1}$.

{\rm (ii)} Since $Y_{i}$ is a nonsingular rational curve, by the adjunction
formula, we have \(K_{S}\cdot Y_{i}+(Y_{i})^{2}=-2\).
Since \(K_{S}=[-Y]\), the conditions (a) and (b) are equivalent to 
each other.

\begin{Lemma}
\label{lem22}
\[H^{2}(F,{\bf Z})\cong\bigoplus_{i=1}^{l} H^{2}(F_{i},{\bf Z})
\cong\bigoplus_{i=1}^{l} {\bf Z}[F_{i}].\]
\end{Lemma}
\noindent
{\em Proof.} Let $F'=\bigcup_{k\neq 1} F_{k}$.
We see that \(H^{1}(F_{1}\cap F',{\bf Z})=
H^{2}(F_{1}\cap F',{\bf Z})=0\).
From the exact sequence (\ref{a3}), we have 
\[H^{2}(F,{\bf Z})\cong H^{2}(F_{1},{\bf Z})\oplus H^{2}(F',{\bf Z}).\]
Using the argument inductively, we obtain the assertion.
\begin{Lemma}
The configuration of $F$ is a tree,
{\rm that is, the dual graph of $F$ is connected and contains no cycles.}
\end{Lemma}
\noindent
{\em Proof.}
Since \(H^{0}(S,{\bf Z})\cong H^{0}(F,{\bf Z})\), the connectivity of $S$
implies the connectivity of $F$.
To show that $F$ contains no cycles,
it is sufficient to show that there are no irreducible components
$F_{i_{1}},\ldots ,F_{i_{m}}\,(m\geq 3)$ such that
\[\begin{array}{ll}
F_{i_{k}}\cap F_{i_{k+1}}\neq \emptyset & (k=1,\ldots ,m-1),\\
F_{i_{m}}\cap F_{i_{1}}\neq \emptyset. &
\end{array}\]
We assume that there exist such $F_{i_{1}},\ldots ,F_{i_{m}}$.
Let \(F^{(1)}=F_{i_{1}}\cup\cdots\cup F_{i_{m}}\) and
\(F^{(2)}=F_{i_{1}}\cup\cdots\cup F_{i_{m-1}}\).
Since $F^{(1)}$ and $F^{(2)}$ are connected, we have $H^{0}(F^{(1)},{\bf Z})
\cong H^{0}(F^{(2)},{\bf Z})\cong{\bf Z}.$
\(F^{(2)}\cap F_{i_{m}}\) consists of at least two points,
and let ${\nu}$  be the number of the points of \(F^{(2)}\cap F_{i_{m}}\).
Then we have \(H^{0}(F^{(2)}\cap F_{i_{m}},{\bf Z})\cong{\bf Z}^{\nu}\,
(\nu\geq 2).\)
Note that \(H^{0}(F_{i_{m}},{\bf Z})\cong H^{0}({\bf P}^{1},{\bf Z})
\cong{\bf Z}.\)
We consider the Mayer--Vietoris exact sequence
\[\begin{array}{cccclcl}
  0 & \rightarrow & \underbrace{H^{0}(F^{(1)},{\bf Z})}_{\cong{\bf Z}}
   & \rightarrow &
  \underbrace{H^{0}(F^{(2)},{\bf Z})\oplus 
  H^{0}(F_{i_{m}},{\bf Z})}_{\cong {\bf Z}\oplus {\bf Z}} 
   & \rightarrow &
  \underbrace{H^{0}(F^{(2)}\cap F_{i_{m}},{\bf Z})}_{\cong {\bf Z}^{\nu}}\\
  & \rightarrow & H^{1}(F^{(1)},{\bf Z}) & \rightarrow 
  & H^{1}(F^{(2)},{\bf Z})\oplus H^{1}(F_{i_{m}},{\bf Z}) & \rightarrow 
  & \ldots.
\end{array}\]
If we assume \(H^{1}(F^{(1)},{\bf Z})=0\) then $\nu$ must be one,
which contradicts the fact that \(\nu\geq 2\).
So we have \(H^{1}(F^{(1)},{\bf Z})\neq 0.\)
Now let $F^{(3)}$ be the union of irreducible components of $F$ which 
do not belong to $F^{(1)}$.
Namely, \(F=F^{(1)}\cup F^{(3)}.\)
We consider the exact sequence
\[\begin{array}{ccccc}
  \underbrace{H^{1}(F,{\bf Z})}_{=0} & \rightarrow & 
  H^{1}(F^{(1)},{\bf Z})\oplus H^{1}(F^{(3)},{\bf Z}) &
  \rightarrow & H^{1}(F^{(1)}\cap F^{(3)},{\bf Z}).
\end{array}\]
Since $F^{(1)}\cap F^{(3)}$ consists of a finite set of points, 
we have \(H^{1}(F^{(1)}\cap F^{(3)},{\bf Z})=0\), 
so \(H^{1}(F^{(1)},{\bf Z})\oplus H^{1}(F^{(3)},{\bf Z})=0.\)
Therefore we have \(H^{1}(F^{(1)},{\bf Z})= 0\),
which is a contradiction.
This completes the proof.\\
\\
{\it Proof of Theorem 2.1.} 
{\rm (i)} Since $S$ is a compactification of ${\bf C}^{2}$, we see that 
$S$ is a rational surface by Theorem 5 of Kodaira \cite{Kod2}.

{\rm (ii)} Let $(S, Y)$ be an Okamoto--Painlev\'{e} pair.  
We know that $S$ is a projective 
rational surface and the configuration of the dual graph of 
$F= S - {\bf C}^2$ is a tree. Let $F = \sum_{i=1}^l F_i$ be the 
irreducible decomposition of $F$, where $l$ denotes the number of the 
irreducible components of $F$. 
Let $[F_{i}]$ denote the class in $H^{2}(S,{\bf Z})$ dual to the class
of $F_{i}\subset S$.
From Lemma~\ref{lem22} and $(\ref{a1})$, we obtain the isomorphism
\[\bigoplus_{i=1}^l {\bf Z}[F_i] \cong H^2(F, {\bf Z}) 
\stackrel{\cong}{\rightarrow}  H^2(S, {\bf Z}).\]
This implies that $\{[F_{i}]\;|\;1\leq i \leq l\}$ forms an integral basis of
$H^{i}(S,{\bf Z})$.

From the  condition (iii) of Definition~\ref{def:op}, we have
$$
Y^2 = \sum_{i=1}^r a_i (Y \cdot Y_i) = 0.
$$
Since $K_S = -Y$, we know that $(K_{S})^{2}=(-Y)^{2}=Y^{2}=0$.
By Noether's formula, we have
$$
e(S) + (K_S)^2 = 12 \chi(S, {\cal O}_S) , 
$$
and $\chi(S, {\cal O}_S) = 1$.
Then we obtain
$$
e(S) = 2 - 2b_1(S) + b_2(S) = 12.
$$
Since $b_1(S) = 0$, we have $b_2(S) =10$, and hence 
$$
\fbox{ $l$ = (the number of the components of $F$) = 10. }
$$
Since $D = \sum_{i=1}^r Y_i \subset F = \sum_{i=1}^{10} F_i$, we have 
$$
r \leq 10.
$$
Since $F = \sum_{i=1}^{10} F_i$ is a connected divisor with normal 
crossings,
the divisor $D = Y_{{\it red}} = \sum_{i=1}^r Y_i$ is also 
a connected divisor with normal crossings. 
To see that $D=Y_{{\it red}}$ is connected, it is sufficient to show that
$Y$ is connected.
By Serre duality, we have 
$H^{1}(S,{\cal O}_{S})
\cong H^{1}(S,{\cal O}(K_{S}))$.
Hence we have $H^{1}(S,{\cal O}(K_{S}))=0$.
By the exact sequence 
\[\begin{array}{ccccccccc}
  0 & \rightarrow & {\cal O}(K_{S}) & \rightarrow & {\cal O}_{S}
  & \rightarrow & {\cal O}_{Y} & \rightarrow & 0,
\end{array}\]
we obtain the long exact sequence
\[\begin{array}{cccclcl}
  0 & \rightarrow & \underbrace{H^{0}(S,{\cal O}(K_{S}))}_{=0}
  & \rightarrow 
  & H^{0}(S,{\cal O}_{S})& \rightarrow & H^{0}(S,{\cal O}_{Y})\\
  & \rightarrow & \underbrace{H^{1}(S,{\cal O}(K_{S}))}_{=0} 
  & \rightarrow & H^{1}(S,{\cal O}_{S}) & \rightarrow & \ldots.
\end{array}\]
Therefore we have 
$H^{0}(S,{\cal O}_{Y})\cong H^{0}(S,{\cal O}_{S})\cong {\bf C}$.
This shows that $Y$ is connected.

Now we prove that the divisor $Y = \sum_{i=1}^r a_i Y_i$ has one of the 
configurations in the list of singular fibers of  elliptic surfaces (cf.\ 
Theorem 6.2, \cite{Kod1}). 
First let us show that the greatest 
common divisor of $\{ a_i \}_{i=1}^r$  is equal to one.
Since $S$ is a rational surface with $b_2(S) = 10$, $S$ is not relatively 
minimal.  Hence $S$ contains an exceptional curve $E$ of the first kind, 
that is, $E \cong {\bf P}^1$ and $E^2 = -1$. 
Note that $K_S = [-Y]$.
By the adjunction formula, we have 
$K_S\cdot E + E^2 = -2$ and this implies that $K_S \cdot E =  -1$ and 
equivalently $Y \cdot E = 1$.  
Since $K_S = [-Y]$ and $Y \cdot Y_i = 0 $ for every irreducible component 
$Y_i$ of $Y$, 
we know that $E$ is not a component of $Y$ hence $E \cdot Y_i 
\geq 0$.
On the other hand, since we have 
$$
1 = Y \cdot E = \sum_{i=1}^r a_i (Y_i \cdot E), 
$$
there exists an irreducible component $Y_i$ with $a_i = 1$.
Under the condition, 
in order to see that the configuration of $Y$ is one in the list of 
the singular fibers of elliptic surfaces in \cite{Kod1},
 we can follow the proof of Theorem 6.2 of \cite{Kod1}.
In fact, only the following conditions are needed to determine all of 
the configuration of singular fibers of elliptic surfaces:
\begin{enumerate}
\item $D = Y_{red}$ is connected.

\item $K_S \cdot Y = 0$.

\item $Y \cdot Y_i = 0$, for $i\;(1 \leq i \leq r)$ .

\end{enumerate}
We have proved the first assertion,
and the third assertion follows from Definition~\ref{def:op}.
The second assertion follows from the third assertion because 
$K_S = [- Y]$.  
Hence the configuration of $Y$ is in the list of Kodaira. 
Next, we note that $D= Y_{red}$ must 
be a divisor with normal crossings and the dual graph of $D$ must be 
a tree.  
Then $Y $  must be one of the types of $\widetilde{E}_{r-1}$ 
for $r =9, 8, 7$ and $\widetilde{D}_{r-1}$ for $r=5,\ldots, 10$.
At this moment, we see that
each irreducible component $Y_i$ of $D = Y_{red}$
is a smooth rational curve with $(Y_i)^2 = -2$. 

We will show that the configurations of types $\widetilde{D}_{9}$ 
and $\widetilde{D}_{8}$  can not occur.  
Let us set $\Lambda_{Y} = \sum_{i=1}^r {\bf Z} [Y_i]$ and 
$\Lambda_{F} = \sum_{j=1}^{10} {\bf Z} [F_j] \cong \Pic (S) \cong 
H^2(S, {\bf Z})$. 
Then we have the inclusion map 
$$
\iota: \Lambda_{Y} \hookrightarrow \Lambda_{F}\cong H^2(S, {\bf Z}).
$$
By this inclusion map $\iota$, $\Lambda_Y$ can be considered as a 
sublattice of $H^2(S, {\bf Z})$ and $H^2(S, {\bf Z})$  is generated by $[F_j] 
\;(1\leq j \leq 10)$. 
The intersection matrix $I_F:= ((F_i\cdot F_j))_{1 \leq i, j \leq 10}$ is 
the unimodular matrix with the 
signature $(b_{+}, b_{-}) = (1, 9)$. On the other hand, the 
intersection matrix
$$
I_Y = ((Y_i\cdot Y_j))_{1 \leq i, j \leq r}
$$
has a null eigenvalue corresponding to $Y$ because $Y^2 = 0$.
This means that 
the rank of $\Lambda_Y$ is strictly less than ten, that is,
$$
r < 10.
$$
This proves that the configuration of type $\widetilde{D}_9$ does not occur. 
If $Y$ is of type $\widetilde{D}_8$, 
then $r = 9$ and hence we can write $F = \sum_{i=1}^9 Y_i + F_{10}$.  
Since  $F$ is connected and the dual graph of $F$ is a tree, we see that 
$F_{10}$ intersects only one irreducible component $Y_i$ for some 
$i\; ( 1 \leq i \leq 9)$. In this case, $\Lambda_F 
= H^2(S, {\bf Z})$ is generated by $Y_1, \ldots, Y_9$ and $F_{10}$ where 
the dual graph of $Y_1, \ldots, Y_9$ is of type $\widetilde{D}_8$. 
Now by a direct calculation we see that the intersection matrix of 
$\{Y_1, \ldots, Y_9, F_{10} \}$ can not have the determinant $-1$.  
Hence $\widetilde{D}_{8}$ does not occur.  
Then we see that the dual graph of $Y$ must be one of the following 
types:
$$
\widetilde{E}_8,\widetilde{E}_7,\widetilde{E}_6, \widetilde{D}_7,
\widetilde{D}_6,\widetilde{D}_5,
\widetilde{D}_4.
$$
Conversely, if $Y$ is one of the types as  above, 
we can construct Okamoto--Painlev\'e pair $(S,Y)$ by blowing up 
and blowing down of $({\bf P}^2, Y = 3H)$ where $H$ denotes a line 
in ${\bf P}^2$. For detail, see \S\ref{s31}. This proves our theorem.

\begin{Remark}\label{remE8}
{\rm 
We will classify not only the configurations of $Y$ 
but also $F$ of an Okamoto--Painlev\'e pair $(S, Y)$. 
By using the result, one can show that $S - D$ is covered by 
a finite number of Zariski open sets $\{ U_i \}$ each of which is  
isomorphic to ${\bf C}^{2}$ if $Y$ is not of type $\widetilde{E}_{8}$.
(See \cite{ST} and \cite{MMT}).
}
\end{Remark} 

\begin{Example}
{\rm 
We consider the case where the configuration of $Y$ is of type 
$\widetilde{D}_4$.
Namely, the corresponding Painlev\'e equation is type $P_{VI}$.
For more details, see \cite{ST}.
At first, we take the Hirzebruch surface ${\Sigma}_{(\varepsilon)}^{(2)}$
which is obtained by gluing four copies of ${\bf C}^2$
via following identification.
Note that ${\Sigma}_{(0)}^{(2)}={\bf F}_{2}$.
}
\[U_{i}=\Spec {{\bf C}[x_{i},y_{i}]} \cong {\bf C}^2\;(i=0,1,2,3)\]
\begin{eqnarray}
x_{0}=x_{1}, && y_{0}=1/y_{1},\nonumber\\
x_{0}=1/x_{2}, && y_{0}=x_{2}(\varepsilon -x_{2}y_{2}),\label{HirSur}\\
x_{2}=x_{3}, && y_{2}=1/y_{3} \nonumber
\end{eqnarray}
\end{Example}
We consider a fiber space $({\Sigma}_{(\varepsilon)}^{(2)}\times B_{VI},
\pi, B_{VI})$,
where $B_{VI}={\bf C}\setminus \{0,1\}$.
Let us take 
\[\varepsilon =(\kappa_{0}+\kappa_{1}+\kappa_{t}-1+\kappa_{\infty})/2,\]
where $\kappa_{\nu}\;(\nu =0,1,t,\infty)$ are complex constants in the 
Hamiltonian function $H_{VI}$
(cf.\ \cite{ST}).
For any parameter $t\in B_{VI}$,
we define a divisor $D^{(0)}(t)$ on ${\Sigma}_{(\varepsilon)}^{(2)}\times t$:
\[D^{(0)}(t)=\{(x_{1},y_{1},t)\in U_{1}\times t \,|\,y_{1}=0\}
\cup \{(x_{3},y_{3},t)\in U_{3}\times t \,|\,y_{3}=0\}.\]
Note that $(D^{(0)}(t))^2 =2$,
and $2D^{(0)}(t)\in |-K_{{\Sigma}_{(\varepsilon)}^{(2)}}|$.
And we take four points 
$a_{\nu}^{(0)}(t)\in D^{(0)}(t)\;(\nu =0,1,t,\infty)$:
\begin{eqnarray*}
a_{\nu}^{(0)}(t) &=& \{(x_{1},y_{1},t)=(\nu,0,t)\}\in U_{1}\cap D^{(0)}(t)
\quad(\nu =0,1,t),\\
a_{\infty}^{(0)}(t) &=& \{(x_{3},y_{3},t)=(0,0,t)\}\in U_{3}\cap D^{(0)}(t).
\end{eqnarray*}
\begin{center}
\begin{picture}(300,150)
\put(60,40){\line(1,0){180}}
\put(60,110){\line(1,0){180}}
\put(75,30){\line(0,1){90}}
\put(225,30){\line(0,1){90}}
\put(75,40){\circle*{5}}
\put(125,40){\circle*{5}}
\put(175,40){\circle*{5}}
\put(225,40){\circle*{5}}
\put(78,25){$a_{0}^{(0)}(t)$}
\put(128,25){$a_{1}^{(0)}(t)$}
\put(178,25){$a_{t}^{(0)}(t)$}
\put(228,25){$a_{\infty}^{(0)}(t)$}
\put(25,35){$D^{(0)}(t)$}
\put(75,110){\vector(1,0){20}}
\put(75,110){\vector(0,-1){20}}
\put(75,40){\vector(1,0){20}}
\put(75,40){\vector(0,1){20}}
\put(225,110){\vector(-1,0){20}}
\put(225,110){\vector(0,-1){20}}
\put(225,40){\vector(-1,0){20}}
\put(225,40){\vector(0,1){20}}
\put(95,100){$x_{0}$}
\put(78,90){$y_{0}$}
\put(95,45){$x_{1}$}
\put(78,60){$y_{1}$}
\put(195,100){$x_{2}$}
\put(212,90){$y_{2}$}
\put(195,45){$x_{3}$}
\put(212,60){$y_{3}$}
\put(120,5){{\bf Figure 2.}}
\end{picture}
\end{center}
\vspace{1cm}

We perform blowings-up to ${\Sigma}_{(\varepsilon)}^{(2)}\times t$
at $a_{\nu}^{(0)}(t)$ for all $t\in B_{VI}$,
and let $D_{\nu}^{(1)}(t)$ be the exceptional curves of the blowings-up
at $a_{\nu}^{(0)}(t)$ for $\nu =0,1,t,\infty$.
We can take four coordinate systems $(z_{\nu},w_{\nu})$ around the points at
infinity of the exceptional curves $D_{\nu}^{(1)}(t)\;(\nu=0,1,t,\infty)$,
where
\begin{eqnarray*}
(z_{\nu},w_{\nu}) & = & ((x_{1}-\nu)y_{1}^{-1},y_{1}) \quad(\nu=0,1,t),\\
(z_{\infty},w_{\infty}) & = & (x_{3}y_{3}^{-1},y_{3}). 
\end{eqnarray*}
Note that we have $w_{\nu}=0$ on $D^{(1)}_{\nu}(t)$ for $\nu=0,1,t,\infty$.
In order to perform the second blowings-up, let us take four points
$a^{(1)}_{\nu}(t)$ for $\nu=0,1,t,\infty$.
\[a_{\nu}^{(1)}(t)=\{(z_{\nu},w_{\nu},t)=(\kappa_{\nu},0,t)\}\in 
D^{(1)}_{\nu}(t)\;(\nu =0,1,t,\infty),\]
\begin{center}
\begin{picture}(280,220)
\thicklines
\put(15,50){\line(1,0){180}}
\put(30,40){\line(0,1){90}}
\put(80,40){\line(0,1){90}}
\put(130,40){\line(0,1){90}}
\put(180,40){\line(0,1){90}}
\thinlines
\multiput(45,180)(5,0){36}{\line(1,0){3}}
\put(25,110){\line(1,2){40}}
\put(175,110){\line(1,2){40}}
\thicklines
\put(60,180){\vector(1,0){20}}
\put(60,180){\vector(-1,-2){10}}
\put(210,180){\vector(-1,0){20}}
\put(210,180){\vector(-1,-2){10}}
\put(30,120){\vector(0,-1){20}}
\put(30,120){\vector(1,2){10}}
\put(80,120){\vector(0,-1){20}}
\put(80,120){\vector(1,2){10}}
\put(130,120){\vector(0,-1){20}}
\put(130,120){\vector(1,2){10}}
\put(180,120){\vector(0,-1){20}}
\put(180,120){\vector(1,2){10}}
\put(60,180){\circle*{4}}
\put(210,180){\circle*{4}}
\put(30,120){\circle*{4}}
\put(80,120){\circle*{4}}
\put(130,120){\circle*{4}}
\put(180,120){\circle*{4}}
\put(65,185){$x_{0}$}
\put(35,165){$y_{0}$}
\put(195,185){$x_{2}$}
\put(210,165){$y_{2}$}
\put(35,110){$z_{0}$}
\put(40,125){$w_{0}$}
\put(85,110){$z_{1}$}
\put(90,125){$w_{1}$}
\put(135,110){$z_{t}$}
\put(140,125){$w_{t}$}
\put(185,110){$z_{\infty}$}
\put(190,125){$w_{\infty}$}
\put(30,80){\circle*{5}}
\put(80,80){\circle*{5}}
\put(130,80){\circle*{5}}
\put(180,80){\circle*{5}}
\put(35,80){$a^{(1)}_{0}(t)$}
\put(85,80){$a^{(1)}_{1}(t)$}
\put(135,80){$a^{(1)}_{t}(t)$}
\put(185,80){$a^{(1)}_{\infty}(t)$}
\put(200,45){$D^{(0)}(t)$}
\put(15,25){$D_{0}^{(1)}(t)$}
\put(65,25){$D_{1}^{(1)}(t)$}
\put(115,25){$D_{t}^{(1)}(t)$}
\put(165,25){$D_{\infty}^{(1)}(t)$}
\put(110,5){{\bf Figure 3.}}
\end{picture}
\end{center}
\vspace{0.5mm}

Let us perform blowings-up at $a_{\nu}^{(1)}(t)$,
and denote $D_{\nu}^{(2)}(t)$ for the exceptional curves, respectively.
We take four coordinate systems $(Z_{\nu},W_{\nu})$ around the points at
infinity of $D_{\nu}^{(2)}(t)$ for $\nu =0,1,t,\infty$, where
\begin{eqnarray*}
(Z_{\nu},W_{\nu}) & = & ((x_1 y_1^{-1}-\kappa_0)y_1^{-1},y_1),\quad 
\nu=0,1,t,\\
(Z_{\infty},W_{\infty}) & = & ((x_3 y_3^{-1}-\kappa_{\infty})y_{3}^{-1},y_3). 
\end{eqnarray*}
\begin{center}
\begin{picture}(400,270)
\thicklines
\put(15,50){\line(1,0){330}}
\put(30,40){\line(0,1){150}}
\put(130,40){\line(0,1){150}}
\put(230,40){\line(0,1){150}}
\put(330,40){\line(0,1){150}}
\thinlines
\multiput(45,240)(5,0){66}{\line(1,0){3}}
\put(25,170){\line(1,2){40}}
\put(325,170){\line(1,2){40}}
\put(20,130){\line(1,0){60}}
\put(120,100){\line(1,0){60}}
\put(220,130){\line(1,0){60}}
\put(320,100){\line(1,0){60}}
\thicklines
\put(60,240){\vector(1,0){20}}
\put(60,240){\vector(-1,-2){10}}
\put(360,240){\vector(-1,0){20}}
\put(360,240){\vector(-1,-2){10}}
\put(80,130){\vector(-1,0){20}}
\put(80,130){\vector(0,1){20}}
\put(180,100){\vector(-1,0){20}}
\put(180,100){\vector(0,1){20}}
\put(280,130){\vector(-1,0){20}}
\put(280,130){\vector(0,1){20}}
\put(380,100){\vector(-1,0){20}}
\put(380,100){\vector(0,1){20}}
\put(60,240){\circle*{4}}
\put(360,240){\circle*{4}}
\put(80,130){\circle*{4}}
\put(180,100){\circle*{4}}
\put(280,130){\circle*{4}}
\put(380,100){\circle*{4}}
\put(70,245){$x_{0}$}
\put(35,225){$y_{0}$}
\put(345,245){$x_{2}$}
\put(360,225){$y_{2}$}
\put(85,135){$W_{0}$}
\put(65,118){$Z_{0}$}
\put(185,105){$W_{1}$}
\put(165,88){$Z_{1}$}
\put(285,135){$W_{t}$}
\put(265,118){$Z_{t}$}
\put(385,105){$W_{\infty}$}
\put(365,88){$Z_{\infty}$}
\put(350,45){$D^{(0)}(t)$}
\put(15,25){$D_{0}^{(1)}(t)$}
\put(115,25){$D_{1}^{(1)}(t)$}
\put(215,25){$D_{t}^{(1)}(t)$}
\put(315,25){$D_{\infty}^{(1)}(t)$}
\put(-10,125){$D_{0}^{(2)}(t)$}
\put(90,95){$D_{1}^{(2)}(t)$}
\put(190,125){$D_{t}^{(2)}(t)$}
\put(290,95){$D_{\infty}^{(2)}(t)$}
\put(170,-5){{\bf Figure 4.}}
\end{picture}
\end{center}
\vspace{0.3cm}

\noindent
For the strict transform of $D_{\nu}^{(i)}(t)$ by the blowing-up,
we also denote by $D_{\nu}^{(i)}(t)$, respectively.
Let $S(t)\rightarrow {\Sigma}_{(\varepsilon)}^{(2)}\times t$ be the 
composition of above eight blowings-up for the parameter $t$.
Then, we see that the configuration of the divisor
\[D(t):=2D^{(0)}(t)+\sum_{\nu=0,1,t,\infty} D_{\nu}^{(1)}(t)\]
on $S(t)$ is of type $\widetilde{D}_4$.
And we see that the complements of $D(t)$ in $S(t)$ is covered by six
Zariski open sets
\begin{eqnarray*}
&& \Spec {\bf C}[Z_{\nu},W_{\nu}] \quad (\nu =0,1,t,\infty),\\
&& \Spec {\bf C}[x_{0},y_{0}],  \\
&& \Spec {\bf C}[x_{2},y_{2}].
\end{eqnarray*}
Note that this example corresponds to the Okamoto--Painlev\'e pair of type
$\widetilde{D}_{4}$--(2) in our classification.
(See \S\ref{s31}).
\begin{Example}
{\rm 
We will construct the space of initial values of the Painlev\'e equation 
of type $P_{I}$
for an example of the Okamoto--Painlev\'e pair of type $\widetilde{E}_{8}$,
As we remarked in Remark \ref{remE8},
the spaces of initial values of $P_{I}$ can not be covered by some Zariski
open sets each of which is isomorphic to ${\bf C}^2$.
For any $t\in B_{I}$, we consider the Hirzebruch surface 
$\Sigma(t):=\Sigma^{(2)}_{(0)}={\bf F}_2$,
and take two curves $D_0(t) = \{y_0=0\}\cup \{y_2=0\}$ 
and $D_0'(t) = \{x_2=0\}\cup \{x_3=0\}$ on $\Sigma(t)$.
Furthermore, we have to consider the multiplicities of each component of the 
anti-canonical divisor.
For the surface $\Sigma(t)$,
we can take the anti-canonical divisor $2D_0(t)+4D_0'(t)\in |-K_{\Sigma(t)}|$.
Let us perform blowing-up at the point $a_0(t)=\{(x_3,y_3)=(0,0)\}$.
\begin{center}
\begin{picture}(300,150)
\put(60,40){\line(1,0){180}}
\put(75,30){\line(0,1){90}}
\put(225,40){\circle*{5}}
\put(228,25){$a_{0}(t)$}
\put(140,115){$D_{0}(t)$}
\put(230,75){$D'_{0}(t)$}
\thicklines
\put(60,110){\line(1,0){180}}
\put(225,30){\line(0,1){90}}
\put(75,110){\vector(1,0){20}}
\put(75,110){\vector(0,-1){20}}
\put(75,40){\vector(1,0){20}}
\put(75,40){\vector(0,1){20}}
\put(225,110){\vector(-1,0){20}}
\put(225,110){\vector(0,-1){20}}
\put(225,40){\vector(-1,0){20}}
\put(225,40){\vector(0,1){20}}
\put(95,100){$x_{0}$}
\put(78,90){$y_{0}$}
\put(95,45){$x_{1}$}
\put(78,60){$y_{1}$}
\put(195,100){$x_{2}$}
\put(212,90){$y_{2}$}
\put(195,45){$x_{3}$}
\put(212,60){$y_{3}$}
\put(215,75){{\bf 4}}
\put(150,100){{\bf 2}}
\put(120,-5){{\bf Figure 5.}}
\end{picture}
\end{center}
\vspace{1.5cm}

\noindent
Let $D_1$ denote the exceptional curve,
and take two coordinate systems $(z_1,W_1)$ and $(Z_1,w_1)$ on $D_1$
which satisfies
$D_1=\{z_1=0\}\cup\{w_1=0\}$ and $Z_1=W_1^{-1}$.
Note that 
\begin{eqnarray*}
(z_1,W_1) & = & (x_3,x_3^{-1} y_3),\\
(Z_1,w_1) & = & (x_3 y_3^{-1},y_3).
\end{eqnarray*}
}
\end{Example}
\begin{center}
\begin{picture}(250,200)
\thinlines
\put(40,70){\line(1,0){110}}
\put(50,60){\line(0,1){120}}
\thicklines
\put(40,170){\line(1,0){160}}
\put(190,110){\line(0,1){70}}
\put(135,65){\line(1,1){60}}
\put(140,70){\vector(1,1){15}}
\put(140,70){\vector(-1,0){20}}
\put(190,120){\vector(0,1){20}}
\put(190,120){\vector(-1,-1){15}}
\put(110,175){$D_0(t)$}
\put(195,140){$D_0'(t)$}
\put(165,85){$D_1(t)$}
\put(118,75){$z_1$}
\put(140,87){$W_1$}
\put(170,115){$Z_1$}
\put(175,130){$w_1$}
\put(120,157){${\bf 2}$}
\put(180,150){${\bf 4}$}
\put(160,100){${\bf 3}$}
\put(100,30){{\bf Figure 6.}}
\end{picture}
\end{center}
Next, let us perform the blowing-up at the point $a_1(t)=\{(Z_1,w_1)=(0,0)\}$,
and denote the exceptional curve by $D_2(t)$.
Take two coordinate systems $(z_2,W_2)$ and $(Z_2,w_2)$ on $D_2$
such that 
$D_2=\{z_2=0\}\cup\{w_2=0\}$ and $Z_2=W_2^{-1}$.
Note that 
\begin{eqnarray*}
(z_2,W_2) & = & (x_3 y_3^{-1},x_3^{-1} y_3^2),\\
(Z_2,w_2) & = & (x_3 y_3^{-2},y_3).
\end{eqnarray*}
\begin{center}
\begin{picture}(250,200)
\thinlines
\put(40,70){\line(1,0){80}}
\put(50,60){\line(0,1){120}}
\thicklines
\put(40,170){\line(1,0){160}}
\put(190,110){\line(0,1){70}}
\put(100,120){\line(1,0){100}}
\put(110,60){\line(0,1){70}}
\put(110,120){\vector(1,0){20}}
\put(110,120){\vector(0,-1){20}}
\put(190,120){\vector(0,1){20}}
\put(190,120){\vector(-1,0){20}}
\put(110,120){\circle*{4}}
\put(190,120){\circle*{4}}
\put(98,127){$Q$}
\put(195,107){$P$}
\put(110,175){$D_0(t)$}
\put(195,140){$D_0'(t)$}
\put(80,85){$D_1(t)$}
\put(70,115){$D_2(t)$}
\put(115,100){$z_2$}
\put(127,107){$W_2$}
\put(165,125){$Z_2$}
\put(175,140){$w_2$}
\put(120,157){${\bf 2}$}
\put(180,150){${\bf 4}$}
\put(115,85){${\bf 3}$}
\put(140,125){${\bf 6}$}
\put(100,30){{\bf Figure 7.}}
\end{picture}
\end{center}
Let $P=P(t)=\{(Z_2,w_2)=(0,0)\}$ and $Q=Q(t)=\{(z_2,W_2)=(0,0)\}$.
Let us perform blowing-up at the point $a_2(t)=\{(Z_2,w_2)=(1/4,0)\}
=\{(z_2,W_2)=(0,4)\}$,
and denote the exceptional curve by $D_3(t)$.
Note that we can consider this blowing-up 
by using the coordinate system
either $(Z_2,w_2)$ or $(z_2,W_2)$.

At first, we consider the blowing-up at $P$ with the coordinate system  $(Z_2,w_2)$.
Take two coordinate systems $(z_3,W_3)$ and $(Z_3,w_3)$ on $D_3(t)$
such that 
$D_3(t)=\{z_3=0\}\cup\{w_3=0\}$ and $Z_3=W_3^{-1}$.
Then we have
\[\left\{\begin{array}{rcl}
    z_3 & = & x_3 y_3^{-2} -1/4,\\
    W_3 & = & (x_3 y_3^{-2} -1/4)^{-1}y_3,
\end{array}\right.\]
and 
\[\left\{\begin{array}{rcl}
    Z_3 & = & (x_3 y_3^{-2} -1/4)y_3^{-1},\\
    w_3 & = & y_3.
\end{array}\right.\]
\vspace{0.8cm}

\begin{center}
\begin{picture}(250,200)
\thinlines
\put(40,70){\line(1,0){80}}
\put(50,60){\line(0,1){120}}
\thicklines
\put(40,170){\line(1,0){160}}
\put(190,110){\line(0,1){70}}
\put(100,120){\line(1,0){100}}
\put(110,60){\line(0,1){70}}
\put(160,50){\line(0,1){80}}
\put(160,120){\vector(-1,0){20}}
\put(160,120){\vector(0,-1){20}}
\put(160,60){\vector(0,1){20}}
\thinlines
\put(160,60){\vector(-1,0){20}}
\put(160,60){\line(1,0){5}}
\put(110,175){$D_0(t)$}
\put(195,140){$D_0'(t)$}
\put(95,135){$D_1(t)$}
\put(70,115){$D_2(t)$}
\put(145,135){$D_3(t)$}
\put(132,112){$z_3$}
\put(142,100){$W_3$}
\put(142,75){$Z_3$}
\put(132,65){$w_3$}
\put(120,157){${\bf 2}$}
\put(180,145){${\bf 4}$}
\put(115,90){${\bf 3}$}
\put(135,125){${\bf 6}$}
\put(165,90){${\bf 5}$}
\put(100,20){{\bf Figure 8.}}
\end{picture}
\end{center}
In the same way, we have to perform five more blowings-up,
and take the coordinate systems $(z_i,W_i)$ and $(Z_i,w_i)$
which satisfies that $D_i(t)=\{z_i=0\}\cup\{w_i=0\}$ and $Z_i=W_i^{-1}$,
where $D_i(t)$ is the exceptional curve of each blowing-up
$(i=4,5,6,7)$.
Similarly, let $a_i(t)$ be the center of each blowing-up.
Then we have to take $a_i(t)$ as follows.
For more details, see \cite{Oka}.
\begin{eqnarray*}
a_3(t) & = & \{(Z_3,w_3)=(0,0)\},\\
a_4(t) & = & \{(Z_4,w_4)=(0,0)\},\\
a_5(t) & = & \{(Z_5,w_5)=(0,0)\},\\
a_6(t) & = & \{(Z_6,w_6)=(t/2,0)\},\\
a_7(t) & = & \{(Z_7,w_7)=(1/2,0)\}.
\end{eqnarray*}
\vspace{6cm}
\begin{center}
\begin{picture}(400,320)
\thinlines
\put(40,190){\line(1,0){80}}
\put(50,180){\line(0,1){120}}
\multiput(270,60)(5,0){16}{\line(1,0){3}}
\thicklines
\put(50,190){\vector(1,0){20}}
\put(50,190){\vector(0,1){20}}
\put(20,205){{\scriptsize $y=y_1$}}
\put(65,182){{\scriptsize $x=x_1$}}
\put(40,290){\line(1,0){160}}
\put(190,230){\line(0,1){70}}
\put(100,240){\line(1,0){100}}
\put(110,180){\line(0,1){70}}
\put(160,170){\line(0,1){80}}
\put(150,180){\line(1,0){80}}
\put(220,110){\line(0,1){80}}
\put(210,120){\line(1,0){80}}
\put(280,50){\line(0,1){80}}
\put(110,240){\vector(1,0){15}}
\put(110,240){\vector(0,-1){15}}
\put(86,225){{\scriptsize $x_3y_3^{-1}$}}
\put(120,245){{\scriptsize $x_3^{-1}y_3^{2}$}}
\put(190,240){\vector(-1,0){15}}
\put(190,240){\vector(0,1){15}}
\put(193,250){{\scriptsize $y_3$}}
\put(165,230){{\scriptsize $x_3y_3^{-2}$}}
\put(110,295){$D_0(t)$}
\put(195,265){$D_0'(t)$}
\put(82,210){$D_1(t)$}
\put(70,238){$D_2(t)$}
\put(165,210){$D_3(t)$}
\put(180,168){$D_4(t)$}
\put(225,145){$D_5(t)$}
\put(240,108){$D_6(t)$}
\put(285,85){$D_7(t)$}
\put(120,277){${\bf 2}$}
\put(180,265){${\bf 4}$}
\put(115,210){${\bf 3}$}
\put(135,230){${\bf 6}$}
\put(150,210){${\bf 5}$}
\put(190,183){${\bf 4}$}
\put(210,145){${\bf 3}$}
\put(250,123){${\bf 2}$}
\put(270,85){${\bf 1}$}
\put(345,60){\vector(-1,0){20}}
\put(345,60){\vector(0,1){20}}
\put(345,60){\circle*{3}}
\put(350,80){$v$}
\put(320,50){$u$}
\put(170,30){{\bf Figure 9.}}
\end{picture}
\end{center}
Let $Y=2D_0+4D'_0+3D_1+6D_2+5D_3+4D_4+3D_5+2D_6+D_7$,
and $D=Y_{\it red}$.
For simplicity, we rewrite $(x_1,y_1)$ as $(x,y)$,
and set $U=\Spec {\bf C}[x,y]$.
Note that 
\[(U;(x,y))\subset S-D.\]
Now we consider the curve $D_8$.
Let us take the coordinate system $(u,v)$ around the point at infinity
of $D_8$.
Let $U'$ be the coordinate system $(U'';(u,v))\cong {\bf C}^2$.
Note that $D_8\cap U'=\{u=0\}$.
By considering the coordinate transformations via above eight blowings-up,
we have the following coordinate transformation between $(U;(x,y))$ and 
$(U';(u,v))$:
\begin{eqnarray}
u & = & x^{15}y^{-8}-\frac{1}{4}x^{12}y^{-6}-\frac{t}{2}x^{4}y^{-2}
+ \frac{1}{2}x^2y^{-1},\label{eq;u}\\
v & = & - x^{-2}y.\label{eq;v}
\end{eqnarray}
By calculating exterior derivations of $u$ and $v$, we have 
\[du\wedge dv=x^{12}y^{-8}dx\wedge dy.\]
On the other hand, by solving the system of equations (\ref{eq;u})
and (\ref{eq;v}), we have 
\begin{eqnarray*}
x & = & A^{-1}v^{-2},\\
y & = & -A^{-2}v^{-3},
\end{eqnarray*}
where
\[A=uv^6+\frac{1}{2}v^5+\frac{t}{2}v^4+\frac{1}{4}.\]
So,
\[x^{12}y^{-8} = A^{4}
=(uv^6+\frac{1}{2}v^5+\frac{t}{2}v^4+\frac{1}{4})^4.\]
Therefore, we have
\[dx\wedge dy=\frac{du\wedge dv}
{(uv^6+\frac{1}{2}v^5+\frac{t}{2}v^4+\frac{1}{4})^4}.\]
\\
\noindent
On the other hand, note that we are denoting $(x_1,y_1)$ by $(x,y)$.
By using $x_1=x_3^{-1}$ and $y_1=-x_3^{-2}y_3$
(see (\ref{HirSur})), 
from (\ref{eq;u}) and (\ref{eq;v}),
we see that 
\begin{eqnarray*}
u & = & x_3 y_3^{-8}-\frac{1}{4}y_3^{-6}-\frac{t}{2}y_3^{-2}
-\frac{1}{2}y_3^{-1},\\
v & = & y_3.
\end{eqnarray*}
Then we have 
\[uv^6 + \frac{1}{2}v^5 +\frac{t}{2}v^4 +\frac{1}{4}
\quad = \quad x_3 y_3^{-2}.\]
From Figure 9, the curve $C=\{uv^6+(1/2)v^5+(t/2)v^4+(1/4)=0\}$ 
coincides with the component $D'_0(t)$ of the divisor $Y$.
Note that the order of the pole at  $C$ of the 2-form $dx_0\wedge dy_0$ 
and the multiplicity of $D'_0(t)$ in $Y$ are both equal to four.
Namely the affine chart $(U';(u,v))\cong {\bf C}^2$ intersects $Y$,
and satisfies that $U'\cap Y=D'_0(t)$.

Next, let us perform the third blowing-up 
(the blowing-up with center of $a_2(t)$), 
by using the coordinate system
whose origin is $Q(t)$, namely $(z_2,W_2)$.
For simplicity, let us use the same notations $(Z_i,w_i)$, $(z_i,W_i)$, 
or $D_i$ as above.
But we denote $(u',v')$ for the coordinate system around 
the point at infinity of $D_8$.
And let $U''$ be the coordinate system $(U'';(u',v'))\cong {\bf C}^2$.
Now we have the following results:
\begin{eqnarray*}
u' & = & x^{-9}y^8 - 4x^{-6} y^6 - \frac{t}{2}x^{-2}y^2 + \frac{1}{2}x^{-1}y\\
   & = & x_3^{-7}y_3^8 - 4x_3^{-6} y_3^6 - \frac{t}{2}x_3^{-2}y_3^2 
         - \frac{1}{2}x_3^{-1}y_3\\
v' & = & -x y^{-1} \\
   & = & -x_3 y_3^{-1}
\end{eqnarray*}
By the similar calculation, we have
\[dx\wedge dy= dx_1\wedge dy_1= 
\frac{du'\wedge dv'}{(u'(v')^6 + \frac{1}{2}(v')^5 
+ \frac{t}{2}(v')^4 + 4 )^3},\]
and
\[u'(v')^6 + \frac{1}{2}(v')^5 + \frac{t}{2}(v')^4 + 4 
\quad = \quad x_3^{-1}y_3^2.\]
From Figure 9, the curve $C'=\{u'(v')^6 + (1/2)(v')^5 + (t/2)(v')^4 + 4 =0\}$ 
coincides with the component $D_1(t)$ of the divisor $Y$.
The order of the pole $C'$ of the 2-form $dx_1\wedge dy_1$ 
and the multiplicity of $D_1(t)$ in $Y$ are both equal to three.
We see that $U''\cap Y=D'_0(t)$.

Consequently, for an Okamoto--Painlev\'e pair of type $\widetilde{E}_8$,
$S-Y$ is covered by $(U;(x_1,y_1))\cong {\bf C}^2$ 
and $(U';(u,v))-C\cong {\bf C}^2-C$ 
(or $(U'';(u',v'))-C'$).
\section{Okamoto--Painlev\'e pairs of non-elliptic type}

\begin{Definition}
{\rm An Okamoto--Painlev\'e pair $(S,Y)$ is 
of {\em elliptic type} if there exists a 
fibration $f:S \longrightarrow {\bf P}^1$ of elliptic curves such that 
a scheme theoretic fiber $f^{*}(\infty)$ at $\infty$ is  $Y$, that is, $f^{*}(\infty) = Y$.  If $(S, Y)$ is not of elliptic type, we 
call $(S, Y)$ is of {\em non-elliptic type}.}
\end{Definition}

A rational elliptic surface can be obtained by blowings-up of 
9 base points of a cubic pencil on ${\bf P}^2$.  
Actually, one can obtain an Okamoto--Painlev\'e pairs $(S, Y)$ 
of type $\tilde{E}_6$ by blowings-up of   infintely near 
base points of a pencil of $ 3 \times line $ and a non-singular 
cubic curve.  This example gives an elliptic fibration 
$f:S \longrightarrow {\bf P}^1$ with $f^{*}(\infty) = Y$, and hence 
$(S, Y)$ is  an Okamoto-Painlev\'e pair of elliptic type.  
On the other hand, if one blows up 
8 (infinitely near) base points of the pencil and  blows up 
a point  on anti-canonical divisor which is not a base points of the pencil, one can not 
obtain an elliptic fibration, so $(S, Y)$ becomes an Okamoto-Painlev\'e pair of non-elliptic type.     
For each type of Okamoto-Painlev\'e pair, one can obtain 
both elliptic type and non-elliptic type depending on 
the position of the points of blowings-up and blowings-down.  
Note that non-elliptic type is general in the 
moduli space of Okamoto-Painlev\'e pairs of each type.

The following porposition is shown in \cite{Sa--T}.

\begin{Proposition}\label{prop/n-e-t}
Let $(S,Y)$ be an Okamoto--Painlev\'e pair.
Then $(S, Y)$ is of non-elliptic type if and only if 
$H^0(S-Y,{\cal O}^{\text{{\it alg}}})\cong {\bf C}$.
(Here $H^0(S-Y,{\cal O}^{\text{{\it alg}}})\cong {\bf C}$ means 
that all of algebraic regular functions on $S-Y$ are
constant functions.)
\end{Proposition}

\section{Construction of type $\widetilde{D}_7$}\label{cstrD7}
Now we will construct the Okamoto--Painlev\'e pair of type $\widetilde{D}_7$
by blowing-up of ${\bf F}_{0}={\bf P}^{1}\times {\bf P}^{1}$ on
$-K_{{\bf F}_{0}}=2s_{0}+2f$.
In the following figure two numbers near each solid line denotes the 
multiplicity and the self-intersection number in $-K_{S}$.
The broken lines denote $(-1)$-curves, whose multiplicities in $F$ are 
zero.
Namely, they are {\em additional components} of $F$.


\vspace{0.8cm}

\begin{center}
\begin{picture}(400,200)(-10, 0)
\thinlines
\put(-20,120){\line(1,0){70}}
\put(-10,110){\line(0,1){70}}
\put(15,109){$(2)$}
\put(18,123){$0$}
\put(-26,147){$(2)$}
\put(-7,147){$0$}
\thicklines
\put(80,130){\vector(-1,0){20}}
\thinlines
\put(120,120){\line(1,0){70}}
\put(130,110){\line(0,1){70}}
\put(80,170){\line(1,0){60}}
\put(180,70){\line(0,1){60}}
\put(150,109){$(2)$}
\put(150,123){$-1$}
\put(114,142){$(2)$}
\put(133,142){$-1$}
\put(100,159){$(1)$}
\put(100,173){$-1$}
\put(164,92){$(1)$}
\put(183,92){$-1$}
\thicklines
\put(220,130){\vector(-1,0){20}}
\thinlines
\put(290,60){\line(0,1){70}}
\put(280,70){\line(1,0){70}}
\put(230,120){\line(1,0){70}}
\put(340,10){\line(0,1){70}}
\put(240,110){\line(0,1){60}}
\put(330,20){\line(1,0){60}}
\put(310,59){$(2)$}
\put(310,73){$-2$}
\put(274,92){$(2)$}
\put(293,92){$-2$}
\put(260,109){$(2)$}
\put(260,123){$-1$}
\put(324,42){$(2)$}
\put(343,42){$-1$}
\put(224,142){$(1)$}
\put(243,142){$-2$}
\put(360,9){$(1)$}
\put(360,23){$-2$}
\end{picture}

\begin{picture}(400,220)(-20,0)
\thicklines
\put(250,200){\vector(0,1){20}}
\thinlines
\put(280,70){\line(0,1){70}}
\put(270,80){\line(1,0){70}}
\put(220,130){\line(1,0){80}}
\put(230,120){\line(0,1){60}}
\put(290,120){\line(0,1){70}}
\put(330,20){\line(0,1){80}}
\put(320,30){\line(1,0){60}}
\put(320,90){\line(1,0){70}}
\put(300,69){$(2)$}
\put(300,83){$-2$}
\put(264,102){$(2)$}
\put(283,102){$-2$}
\put(250,119){$(2)$}
\put(250,133){$-2$}
\put(314,52){$(2)$}
\put(333,52){$-2$}
\put(214,152){$(1)$}
\put(233,152){$-2$}
\put(350,19){$(1)$}
\put(350,33){$-2$}
\put(274,152){$(1)$}
\put(293,152){$-1$}
\put(350,79){$(1)$}
\put(350,93){$-1$}
\thicklines
\put(190,100){\vector(1,0){20}}
\thinlines
\put(60,70){\line(0,1){70}}
\put(50,80){\line(1,0){70}}
\put(0,130){\line(1,0){80}}
\put(10,120){\line(0,1){60}}
\put(70,120){\line(0,1){70}}
\put(110,20){\line(0,1){80}}
\put(100,30){\line(1,0){60}}
\put(100,90){\line(1,0){70}}

\put(160,80){\line(0,1){3}}
\put(160,85){\line(0,1){3}}
\put(160,90){\line(0,1){3}}
\put(160,95){\line(0,1){3}}
\put(160,100){\line(0,1){3}}
\put(160,105){\line(0,1){3}}
\put(160,110){\line(0,1){3}}
\put(160,115){\line(0,1){3}}
\put(160,120){\line(0,1){3}}
\put(160,125){\line(0,1){3}}
\put(160,130){\line(0,1){3}}
\put(160,135){\line(0,1){3}}

\put(60,180){\line(1,0){3}}
\put(65,180){\line(1,0){3}}
\put(70,180){\line(1,0){3}}
\put(75,180){\line(1,0){3}}
\put(80,180){\line(1,0){3}}
\put(85,180){\line(1,0){3}}
\put(90,180){\line(1,0){3}}
\put(95,180){\line(1,0){3}}
\put(100,180){\line(1,0){3}}
\put(105,180){\line(1,0){3}}
\put(110,180){\line(1,0){3}}
\put(115,180){\line(1,0){3}}
\put(80,69){$(2)$}
\put(80,83){$-2$}
\put(44,102){$(2)$}
\put(63,102){$-2$}
\put(30,119){$(2)$}
\put(30,133){$-2$}
\put(94,52){$(2)$}
\put(113,52){$-2$}
\put(-6,152){$(1)$}
\put(13,152){$-2$}
\put(130,19){$(1)$}
\put(130,33){$-2$}
\put(54,152){$(1)$}
\put(73,152){$-2$}
\put(130,79){$(1)$}
\put(130,93){$-2$}
\put(163,112){$-1$}
\put(90,183){$-1$}
\put(190,0){{\bf Figure 10.}}
\end{picture}
\end{center}


\vspace{0.4cm}


\section{Additional Components}
\label{s31}

By definition, 
each Okamoto--Painlev\'e pair $(S, Y)$  
contains the affine plain ${\bf C}^2$ as a Zariski open 
set and we set 
$$
F = S - {\bf C}^2.
$$
We can see $F$ as a divisor with 
ten components which is obtained by
adding some irreducible smooth rational curves to $Y_{red}$.
We call such curves {\em additional components}.

Now we will consider $F$ counting the multiplicity 
$m_{j}$ for the component   
$F_{j}\; (j=1,\ldots ,10)$.
We set $F=\sum_{j=1}^{10} m_{j} F_{j}$ and $Y=\sum_{i=1}^{r} a_{i} Y_{i}$.
Note that $r<10$.
If $F_{j}=Y_{i}$, then we set $m_{j}:=a_{i}$.
And if $F_{j}$ is an additional component, then we set $m_{j}=0$.
\begin{Lemma}
\label{lem31}
Let $C$ be an additional component of $F$ which intersects $F_{j}$.
Then we have \(C^{2}=m_{j}-2\)
\end{Lemma}
\noindent
{\em Proof.} Note that $F_{j}$ which satisfies $m_{j}>0$ and which $C$ 
intersects exists uniquely.
By the adjunction formula, we have $C^{2}=C\cdot Y-2$.
Since $C\cdot Y=m_{j}$, our assertion is proved.\\
\\
Now we assume that $F=\sum_{i=1}^{10} F_i$ satisfies the following condition.\\
\\
$(\ast)$ Both $(F_{1})^{2}$ and $(F_{2})^{2}$ are even numbers, and
both $F_{1}$ and $F_{2}$ intersect only $F_{3}$.

\begin{center}
\begin{picture}(150,200)
\put(80,160){\line(1,0){50}}
\put(80,90){\line(1,0){50}}
\put(90,80){\line(0,1){90}}
\put(100,125){\line(-1,0){50}}
\put(45,125){\circle*{3}}
\put(40,125){\circle*{3}}
\put(35,125){\circle*{3}}
\put(30,125){\circle*{3}}
\put(135,155){$F_{1}$}
\put(135,85){$F_{2}$}
\put(85,65){$F_{3}$}
\put(45,30){{\bf Figure 11.}}
\end{picture}
\end{center}

\begin{Lemma}
\label{lem32}
There exists no Okamoto--Painlev\'e pairs which satisfy $(\ast)$.
\end{Lemma}
\noindent
{\em Proof.} 
\vspace{0.5cm}
We consider the determinant of the intersection matrix of $F$.

\[\begin{array}{ccc}
\det I_{F} & = &
\left| \begin{array}{ccccccc}
(F_{1})^{2} & 0 & 1 & 0 & \cdots & 0 \\
0 & (F_{2})^{2} & 1 & 0 & \cdots & 0 \\
1 & 1 & (F_{3})^{2} & \ast & \cdots & \ast \\
0 & 0 & \ast & \ast & \cdots & \ast \\
\vdots & \vdots & \vdots & \vdots& \ddots & \vdots \\
0 & 0 & \ast & \ast & \cdots & \ast \end{array}
\right| 
\end{array}\]
\vspace{0.5cm}

We see that $\det I_{F}$ must be an even number, if we expand the determinant.
Namely $\det I_{F}\neq -1$.
This contradicts the fact that $I_{F}$ is a unimodular matrix with the 
signature $(b_{+},b_{-})=(1,9)$.
This completes the proof.
\begin{Lemma}
There is not more than one additional component which intersects
the components of $F$ whose multiplicities are $\geq 2$.
\end{Lemma}
\noindent
{\em Proof.} We consider additional components $C_{1}$ and $C_{2}$.
We assume that $C_{1}$ and $C_{2}$ intersect $F_{1}$ and $F_{2}$ respectively,
and that $m_{1}\geq 2$ and $m_{2}\geq 2$.
Note that $C_{1}\cdot C_{2}=0$.
If $C_{1}\cdot C_{2}\neq 0$, then $F$ includes a cycle.
From Lemma~\ref{lem31}, it follows that $(C_{1})^{2}=m_{1}-2$ and 
$(C_{2})^{2}=m_{2}-2$.
Note that $(C_{1})^{2}\geq 0$ and $(C_{2})^{2}\geq 0$.
We will prove by separating into some cases.

{\rm [i]} {\em The case where \(m_{1}\geq 3\) or \(m_{2}\geq 3\).}\\
For example, we assume $m_{1}\geq 3$.
We have $(C_{1})^{2}>0$ and $C_{1}\cdot C_{2}=0$.
Therefore, the Hodge index theorem shows that $(C_{2})^{2}<0$
or $C_{2}\equiv 0$.
This is a contradiction.

{\rm [ii]} {\em The case where $m_{1}=m_{2}=2$.}

{\rm (1)} {\em The case where another additional component intersects 
$C_{1}$ or $C_{2}$.}\\
In this case, we suppose that $C_{3}$ intersects $C_{1}$ as an example.
Then we have $(C_{3})^{2}=-2, \,C_{1}\cdot C_{3}=1$ and 
$C_{2}\cdot C_{3}=0$.
Now we consider a divisor $A=nC_{1}+C_{3}$ for $n\in {\bf Z}_{\geq 2}$.
We have $A^{2}=2n-2>0$ and $A\cdot C_{2}=0$.
By applying the Hodge index theorem again, it follows that $(C_{2})^{2}<0$
or $C_{2}\equiv 0$.
Thus we obtained a contradiction.

{\rm (2)} {\em The case where another additional component intersect neither 
$C_{1}$ nor $C_{2}$.}\\
In this case, $C_{1}$ and $C_{2}$ can not intersect the same component of $F$
by Lemma \ref{lem32}.
Hence we see that $C_{1}$ and $C_{2}$ intersects the different components.
Let us denote them by $F_{1}$ and $F_{2}$ respectively.
Note that $(C_{1})^{2}=(C_{2})^{2}=0,\,C_{1}\cdot F_{1}=1,\,
C_{2}\cdot F_{1}=0$ and $(F_{1})^{2}=-2$.
Now we consider a divisor $B=nC_{1}+F_{1}$ for $n\in {\bf Z}_{\geq 2}$.
Then we have $B^{2}=2n-2>0$ and $B\cdot C_{2}=0$.
This contradicts the Hodge index theorem.
This completes the proof.

\vspace{1cm}


By using the above lemmas, we will be able to give  
the complete list of configurations of $F$ and $Y \subset F$. 

Moreover, one can also check that  
each pattern in the following list can be transformed into 
$({\bf P}^2, -K_{{\bf P}^2}=3h)$,
$({\bf F}_0, -K_{{\bf F}_0}=2s_{0}+2f)$ 
or $({\bf F}_2, -K_{{\bf F}_2}=2s_{\infty})$
by performing blowing-up and blowing-down to $S$ on $-K_{S}$. 
Since $({\bf F}_0, -K_{{\bf F}_0}=2s_{0}+2f)$ 
or $({\bf F}_2, -K_{{\bf F}_2}=2s_{\infty})$
can be transformed into $({\bf P}^2, -K_{{\bf P}^2}=3h)$
by  birational transformations, 
as a consequence, each pattern can be transformed into 
$({\bf P}^2, -K_{{\bf P}^2}=3h)$ by 
performing blowing-up and blowing-down of $(S, Y)$. 
(This gives a proof of the assertion (iii) of Theorem \ref{clop}.) 
This also shows that all of the configuration of $F$ in the 
following list do exist.  In fact, one can perform 
the inverse birationl transformation from $({\bf P}^2, -3 l)$ 
to $(S, Y)$ or $(S, F)$.

\vspace{1cm}

\begin{center}
\large
{\bf COMPLETE LIST OF CONFIGURATIONS OF $F$.}
\end{center}
\normalsize
We will give the list of configurations of $F$ in the form of dual graphs.
In each figure, we denote by each vertex a component with the positive 
multiplicity in $F$, and by each circle an additional component.
The number near each vertex means the multiplicity of the component in $F$,
and the number inside each circle means the self-intersection number of the corresponding 
additional component.
Note that the multiplicity of an additional component in $F$ is zero.


\begin{picture}(350,170)(-50, 0)
\put(-50,150){\fbox{ Type $\widetilde{D}_4$ }}
\put(20,80){\line(1,1){38}}
\put(20,80){\line(-1,1){38}}
\put(20,80){\line(-1,-1){38}}
\put(20,80){\line(1,-1){38}}
\put(20,80){\line(-1,0){45}}
\thicklines
\put(65,125){\circle{20}}
\put(-25,125){\circle{20}}
\put(-25,35){\circle{20}}
\put(65,35){\circle{20}}
\put(-35,80){\circle{20}}

\put(45,105){\circle*{5}}
\put(-5,105){\circle*{5}}
\put(-5,55){\circle*{5}}
\put(45,55){\circle*{5}}
\put(20,80){\circle*{5}}

\put(10,0){$\widetilde{D}_4$--(1)}

\put(58,122){$-1$}
\put(-32,122){$-1$}
\put(-32,32){$-1$}
\put(58,32){$-1$}
\put(-38,77){$0$}

\thinlines
\put(250,80){\line(1,1){25}}
\put(250,80){\line(-1,1){38}}
\put(250,80){\line(-1,-1){38}}
\put(250,80){\line(1,-1){38}}
\put(275,105){\line(1,0){15}}
\put(275,105){\line(0,1){15}}
\thicklines
\put(295,35){\circle{20}}
\put(205,125){\circle{20}}
\put(205,35){\circle{20}}
\put(300,105){\circle{20}}
\put(275,130){\circle{20}}
\put(275,105){\circle*{5}}
\put(225,105){\circle*{5}}
\put(225,55){\circle*{5}}
\put(275,55){\circle*{5}}
\put(250,80){\circle*{5}}

\put(240,0){$\widetilde{D}_4$--(2)}

\put(288,32){$-1$}
\put(198,122){$-1$}
\put(198,32){$-1$}
\put(268,127){$-1$}
\put(293,102){$-1$}

\put(36,105){1}
\put(-1,105){1}
\put(-1,49){1}
\put(36,49){1}
\put(18,85){2}

\put(266,105){1}
\put(229,105){1}
\put(229,49){1}
\put(266,49){1}
\put(248,85){2}

\end{picture}

\vspace{1.5cm}

\begin{picture}(400,170)(-50, 0)

\thinlines
\put(20,80){\line(1,1){38}}
\put(20,80){\line(-1,1){38}}
\put(20,80){\line(-1,-1){38}}
\put(20,80){\line(1,-1){38}}
\put(72,132){\line(1,1){11}}
\thicklines
\put(65,125){\circle{20}}
\put(-25,125){\circle{20}}
\put(-25,35){\circle{20}}
\put(65,35){\circle{20}}
\put(90,150){\circle{20}}

\put(45,105){\circle*{5}}
\put(-5,105){\circle*{5}}
\put(-5,55){\circle*{5}}
\put(45,55){\circle*{5}}
\put(20,80){\circle*{5}}

\put(10,0){$\widetilde{D}_4$--(3)}

\put(58,122){$-1$}
\put(-32,122){$-1$}
\put(-32,32){$-1$}
\put(58,32){$-1$}
\put(83,147){$-2$}

\thinlines
\put(300,80){\line(-1,1){38}}
\put(300,80){\line(-1,-1){38}}
\put(300,80){\line(1,-1){25}}
\put(300,80){\line(1,1){25}}
\put(300,80){\line(-1,0){45}}
\put(235,80){\line(-1,0){10}}
\put(325,105){\line(1,1){13}}
\thicklines
\put(255,125){\circle{20}}
\put(255,35){\circle{20}}
\put(345,125){\circle{20}}
\put(245,80){\circle{20}}
\put(215,80){\circle{20}}

\put(325,105){\circle*{5}}
\put(275,105){\circle*{5}}
\put(275,55){\circle*{5}}
\put(325,55){\circle*{5}}
\put(300,80){\circle*{5}}

\put(290,0){$\widetilde{D}_4$--(4)}

\put(248,122){$-1$}
\put(248,32){$-1$}
\put(338,122){$-1$}
\put(242,77){$0$}
\put(208,77){$-2$}

\put(36,105){1}
\put(-1,105){1}
\put(-1,49){1}
\put(36,49){1}
\put(18,85){2}

\put(316,105){1}
\put(279,105){1}
\put(279,49){1}
\put(316,49){1}
\put(298,85){2}

\end{picture}

\vspace{1cm}

\begin{picture}(400,150)(-50, 0)
\thinlines
\put(20,80){\line(-1,1){38}}
\put(20,80){\line(-1,-1){38}}
\put(20,80){\line(1,-1){25}}
\put(20,80){\line(1,1){43}}
\put(45,105){\line(0,1){20}}
\put(45,105){\line(1,0){20}}
\thicklines
\put(-25,125){\circle{20}}
\put(-25,35){\circle{20}}
\put(70,130){\circle{20}}
\put(45,135){\circle{20}}
\put(75,105){\circle{20}}

\put(45,105){\circle*{5}}
\put(-5,105){\circle*{5}}
\put(-5,55){\circle*{5}}
\put(45,55){\circle*{5}}
\put(20,80){\circle*{5}}

\put(10,0){$\widetilde{D}_4$--(5)}

\put(-32,122){$-1$}
\put(-32,32){$-1$}
\put(63,127){$-1$}
\put(38,132){$-1$}
\put(68,102){$-1$}

\thinlines
\put(300,80){\line(-1,1){38}}
\put(300,80){\line(-1,-1){38}}
\put(300,80){\line(1,-1){25}}
\put(300,80){\line(1,1){25}}
\put(325,105){\line(0,1){15}}
\put(325,105){\line(1,0){15}}
\put(360,105){\line(1,0){10}}
\thicklines
\put(255,125){\circle{20}}
\put(255,35){\circle{20}}
\put(380,105){\circle{20}}
\put(325,130){\circle{20}}
\put(350,105){\circle{20}}

\put(325,105){\circle*{5}}
\put(275,105){\circle*{5}}
\put(275,55){\circle*{5}}
\put(325,55){\circle*{5}}
\put(300,80){\circle*{5}}

\put(290,0){$\widetilde{D}_4$--(6)}

\put(248,122){$-1$}
\put(248,32){$-1$}
\put(373,102){$-2$}
\put(318,127){$-1$}
\put(343,102){$-1$}

\put(36,105){1}
\put(-1,105){1}
\put(-1,49){1}
\put(36,49){1}
\put(18,85){2}

\put(316,105){1}
\put(279,105){1}
\put(279,49){1}
\put(316,49){1}
\put(298,85){2}

\end{picture}

\vspace{1cm}

\begin{picture}(400,200)(-50,0)

\thinlines
\put(20,80){\line(1,1){38}}
\put(20,80){\line(-1,1){38}}
\put(20,80){\line(-1,-1){38}}
\put(20,80){\line(1,-1){25}}
\put(72,132){\line(1,1){11}}
\put(97,157){\line(1,1){11}}
\thicklines
\put(65,125){\circle{20}}
\put(-25,125){\circle{20}}
\put(-25,35){\circle{20}}
\put(115,175){\circle{20}}
\put(90,150){\circle{20}}

\put(45,105){\circle*{5}}
\put(-5,105){\circle*{5}}
\put(-5,55){\circle*{5}}
\put(45,55){\circle*{5}}
\put(20,80){\circle*{5}}

\put(10,0){$\widetilde{D}_4$--(7)}

\put(58,122){$-1$}
\put(-32,122){$-1$}
\put(-32,32){$-1$}
\put(108,172){$-2$}
\put(83,147){$-2$}

\put(36,105){1}
\put(-1,105){1}
\put(-1,49){1}
\put(36,49){1}
\put(18,85){2}

\end{picture}

\vspace{1cm}

\begin{picture}(400,170)(-50,0)

\put(-50,150){\fbox{ Type $\widetilde{D}_5$ }}

\thinlines
\put(45,80){\line(1,1){38}}
\put(20,80){\line(-1,1){38}}
\put(20,80){\line(-1,-1){38}}
\put(45,80){\line(1,-1){38}}
\put(20,80){\line(1,0){25}}
\thicklines
\put(90,125){\circle{20}}
\put(-25,125){\circle{20}}
\put(-25,35){\circle{20}}
\put(90,35){\circle{20}}

\put(70,105){\circle*{5}}
\put(-5,105){\circle*{5}}
\put(-5,55){\circle*{5}}
\put(70,55){\circle*{5}}
\put(20,80){\circle*{5}}
\put(45,80){\circle*{5}}

\put(25,0){$\widetilde{D}_5$--(1)}

\put(83,122){$-1$}
\put(-32,122){$-1$}
\put(-32,32){$-1$}
\put(83,32){$-1$}
\put(-1,105){1}
\put(-1,49){1}
\put(61,105){1}
\put(61,49){1}
\put(18,85){2}
\put(43,85){2}

\thinlines
\put(255,80){\line(-1,1){38}}
\put(255,80){\line(-1,-1){25}}
\put(280,80){\line(1,-1){38}}
\put(280,80){\line(1,1){25}}
\put(305,105){\line(0,1){15}}
\put(305,105){\line(1,0){15}}
\put(280,80){\line(-1,0){25}}
\thicklines
\put(210,125){\circle{20}}
\put(325,35){\circle{20}}
\put(305,130){\circle{20}}
\put(330,105){\circle{20}}

\put(305,105){\circle*{5}}
\put(230,105){\circle*{5}}
\put(230,55){\circle*{5}}
\put(305,55){\circle*{5}}
\put(280,80){\circle*{5}}
\put(255,80){\circle*{5}}
\put(260,0){$\widetilde{D}_5$--(2)}

\put(203,122){$-1$}
\put(318,32){$-1$}
\put(298,127){$-1$}
\put(323,102){$-1$}
\put(234,105){1}
\put(234,49){1}
\put(296,105){1}
\put(296,49){1}
\put(253,85){2}
\put(278,85){2}
\end{picture}

\vspace{1cm}

\begin{picture}(400,170)(-50,0)
\thinlines
\put(45,80){\line(1,1){38}}
\put(20,80){\line(-1,1){38}}
\put(20,80){\line(-1,-1){25}}
\put(45,80){\line(1,-1){38}}
\put(20,80){\line(1,0){25}}
\put(97,132){\line(1,1){11}}
\thicklines
\put(90,125){\circle{20}}
\put(-25,125){\circle{20}}
\put(115,150){\circle{20}}
\put(90,35){\circle{20}}

\put(70,105){\circle*{5}}
\put(-5,105){\circle*{5}}
\put(-5,55){\circle*{5}}
\put(70,55){\circle*{5}}
\put(20,80){\circle*{5}}
\put(45,80){\circle*{5}}

\put(25,0){$\widetilde{D}_5$--(3)}

\put(83,122){$-1$}
\put(-32,122){$-1$}
\put(108,147){$-2$}
\put(83,32){$-1$}
\put(-1,105){1}
\put(-1,49){1}
\put(61,105){1}
\put(61,49){1}
\put(18,85){2}
\put(43,85){2}

\thinlines
\put(280,80){\line(1,1){38}}
\put(255,80){\line(-1,1){38}}
\put(255,80){\line(-1,-1){25}}
\put(280,80){\line(1,-1){25}}
\put(255,80){\line(1,0){25}}
\put(255,80){\line(0,1){35}}
\put(230,55){\line(-1,-1){13}}
\thicklines
\put(325,125){\circle{20}}
\put(210,125){\circle{20}}
\put(255,125){\circle{20}}
\put(210,35){\circle{20}}

\put(305,105){\circle*{5}}
\put(230,105){\circle*{5}}
\put(230,55){\circle*{5}}
\put(305,55){\circle*{5}}
\put(255,80){\circle*{5}}
\put(280,80){\circle*{5}}

\put(260,0){$\widetilde{D}_5$--(4)}

\put(318,122){$-1$}
\put(203,122){$-1$}
\put(252,122){$0$}
\put(203,32){$-1$}
\put(234,105){1}
\put(234,49){1}
\put(296,105){1}
\put(296,49){1}
\put(253,65){2}
\put(278,65){2}

\end{picture}

\vspace{1cm}

\begin{picture}(400,170)(-50,0)

\thinlines
\put(20,80){\line(-1,1){38}}
\put(20,80){\line(-1,-1){25}}
\put(45,80){\line(1,-1){25}}
\put(45,80){\line(1,1){25}}
\put(45,80){\line(-1,0){25}}
\put(20,80){\line(0,1){35}}
\put(20,135){\line(0,1){10}}
\put(70,105){\line(1,1){13}}
\thicklines
\put(-25,125){\circle{20}}
\put(90,125){\circle{20}}
\put(20,125){\circle{20}}
\put(20,155){\circle{20}}

\put(70,105){\circle*{5}}
\put(-5,105){\circle*{5}}
\put(-5,55){\circle*{5}}
\put(70,55){\circle*{5}}
\put(45,80){\circle*{5}}
\put(20,80){\circle*{5}}
\put(25,0){$\widetilde{D}_5$--(5)}

\put(-32,122){$-1$}
\put(83,122){$-1$}
\put(17,122){$0$}
\put(13,152){$-2$}

\put(-1,105){1}
\put(-1,49){1}
\put(61,105){1}
\put(61,49){1}
\put(18,65){2}
\put(43,65){2}

\thinlines
\put(255,80){\line(-1,1){38}}
\put(255,80){\line(-1,-1){25}}
\put(280,80){\line(1,-1){25}}
\put(280,80){\line(1,1){25}}
\put(305,105){\line(0,1){20}}
\put(305,105){\line(1,0){20}}
\put(280,80){\line(-1,0){25}}
\put(305,105){\line(1,1){18}}
\thicklines
\put(210,125){\circle{20}}
\put(305,135){\circle{20}}
\put(335,105){\circle{20}}
\put(330,130){\circle{20}}

\put(305,105){\circle*{5}}
\put(230,105){\circle*{5}}
\put(230,55){\circle*{5}}
\put(305,55){\circle*{5}}
\put(280,80){\circle*{5}}
\put(255,80){\circle*{5}}

\put(260,0){$\widetilde{D}_5$--(6)}

\put(203,122){$-1$}
\put(298,132){$-1$}
\put(328,102){$-1$}
\put(323,127){$-1$}

\put(234,105){1}
\put(234,49){1}
\put(296,105){1}
\put(296,49){1}
\put(253,85){2}
\put(278,85){2}

\end{picture}

\vspace{1cm}

\begin{picture}(400,200)(-50,0)
\thinlines
\put(20,80){\line(-1,1){38}}
\put(20,80){\line(-1,-1){25}}
\put(45,80){\line(1,-1){25}}
\put(45,80){\line(1,1){25}}
\put(20,80){\line(1,0){25}}
\put(70,105){\line(0,1){15}}
\put(70,105){\line(1,0){15}}
\put(105,105){\line(1,0){10}}

\thicklines
\put(-25,125){\circle{20}}
\put(70,130){\circle{20}}
\put(95,105){\circle{20}}
\put(125,105){\circle{20}}

\put(20,80){\circle*{5}}
\put(45,80){\circle*{5}}
\put(-5,105){\circle*{5}}
\put(-5,55){\circle*{5}}
\put(70,105){\circle*{5}}
\put(70,55){\circle*{5}}

\put(25,0){$\widetilde{D}_5$--(7)}

\put(-32,122){$-1$}
\put(63,127){$-1$}
\put(88,102){$-1$}
\put(118,102){$-2$}

\put(-1,105){1}
\put(-1,49){1}
\put(61,105){1}
\put(61,49){1}
\put(18,85){2}
\put(43,85){2}

\thinlines
\put(255,80){\line(-1,1){38}}
\put(255,80){\line(-1,-1){25}}
\put(280,80){\line(1,-1){25}}
\put(280,80){\line(1,1){25}}
\put(255,80){\line(1,0){25}}
\put(305,105){\line(1,1){13}}
\put(332,132){\line(1,1){11}}
\put(357,157){\line(1,1){11}}
\thicklines
\put(210,125){\circle{20}}
\put(325,125){\circle{20}}
\put(350,150){\circle{20}}
\put(375,175){\circle{20}}

\put(305,105){\circle*{5}}
\put(230,105){\circle*{5}}
\put(230,55){\circle*{5}}
\put(305,55){\circle*{5}}
\put(280,80){\circle*{5}}
\put(255,80){\circle*{5}}
\put(260,0){$\widetilde{D}_5$--(8)}

\put(203,122){$-1$}
\put(318,122){$-1$}
\put(343,147){$-2$}
\put(368,172){$-2$}

\put(234,105){1}
\put(234,49){1}
\put(296,105){1}
\put(296,49){1}
\put(253,85){2}
\put(278,85){2}

\end{picture}

\vspace{1.5cm}

\begin{picture}(400,170)(-50,0)
\put(-50,150){\fbox{ Type $\widetilde{D}_6$ }}
\thinlines

\put(20,80){\line(-1,1){38}}
\put(20,80){\line(-1,-1){25}}
\put(70,80){\line(1,-1){25}}
\put(70,80){\line(1,1){25}}
\put(20,80){\line(1,0){50}}
\put(95,105){\line(1,1){13}}
\put(-5,55){\line(-1,-1){13}}

\thicklines
\put(-25,125){\circle{20}}
\put(115,125){\circle{20}}
\put(-25,35){\circle{20}}

\put(95,105){\circle*{5}}
\put(-5,105){\circle*{5}}
\put(-5,55){\circle*{5}}
\put(95,55){\circle*{5}}
\put(70,80){\circle*{5}}
\put(45,80){\circle*{5}}
\put(20,80){\circle*{5}}
\put(30,0){$\widetilde{D}_6$--(1)}

\put(-32,122){$-1$}
\put(108,122){$-1$}
\put(-32,32){$-1$}

\thinlines
\put(250,80){\line(-1,1){38}}
\put(250,80){\line(-1,-1){25}}
\put(300,80){\line(1,-1){25}}
\put(300,80){\line(1,1){25}}
\put(250,80){\line(1,0){50}}
\put(275,80){\line(0,1){35}}
\put(325,105){\line(1,1){13}}

\thicklines
\put(205,125){\circle{20}}
\put(345,125){\circle{20}}
\put(275,125){\circle{20}}

\put(325,105){\circle*{5}}
\put(225,105){\circle*{5}}
\put(225,55){\circle*{5}}
\put(325,55){\circle*{5}}
\put(300,80){\circle*{5}}
\put(275,80){\circle*{5}}
\put(250,80){\circle*{5}}
\put(270,0){$\widetilde{D}_6$--(2)}

\put(198,122){$-1$}
\put(338,122){$-1$}
\put(272,122){$0$}

\put(-1,105){1}
\put(-1,49){1}
\put(86,105){1}
\put(86,49){1}
\put(18,65){2}
\put(43,65){2}
\put(68,65){2}
\put(229,105){1}
\put(229,49){1}
\put(316,105){1}
\put(316,49){1}
\put(248,65){2}
\put(273,65){2}
\put(298,65){2}

\end{picture}

\vspace{2cm}

\begin{picture}(400,170)(-50,0)

\put(-50,160){\fbox{ Type $\widetilde{D}_7$ }}

\thinlines
\put(20,80){\line(-1,1){38}}
\put(20,80){\line(-1,-1){25}}
\put(95,80){\line(1,-1){25}}
\put(95,80){\line(1,1){25}}
\put(20,80){\line(1,0){75}}
\put(120,105){\line(1,1){13}}

\thicklines
\put(-25,125){\circle{20}}
\put(140,125){\circle{20}}

\put(120,105){\circle*{5}}
\put(-5,105){\circle*{5}}
\put(-5,55){\circle*{5}}
\put(120,55){\circle*{5}}
\put(95,80){\circle*{5}}
\put(70,80){\circle*{5}}
\put(45,80){\circle*{5}}
\put(20,80){\circle*{5}}

\put(-32,122){$-1$}
\put(133,122){$-1$}

\put(-1,105){1}
\put(-1,49){1}
\put(111,105){1}
\put(111,49){1}
\put(18,65){2}
\put(43,65){2}
\put(68,65){2}
\put(93,65){2}

\end{picture}

\vspace{1cm}

\begin{picture}(400,150)(-50,0)


\put(-50,140){\fbox{ Type $\widetilde{E}_6$ }}

\thinlines
\put(-30,40){\line(1,0){100}}
\put(20,40){\line(0,1){50}}
\put(-30,40){\line(-1,0){15}}
\put(70,40){\line(1,0){15}}
\put(20,90){\line(0,1){15}}

\thicklines
\put(-55,40){\circle{20}}
\put(95,40){\circle{20}}
\put(20,115){\circle{20}}

\put(-30,40){\circle*{5}}
\put(-5,40){\circle*{5}}
\put(20,40){\circle*{5}}
\put(45,40){\circle*{5}}
\put(70,40){\circle*{5}}
\put(20,65){\circle*{5}}
\put(20,90){\circle*{5}}

\put(-62,37){$-1$}
\put(88,37){$-1$}
\put(13,112){$-1$}

\put(10,0){$\widetilde{E}_6$--(1)}

\thinlines
\put(250,40){\line(1,0){100}}
\put(300,40){\line(0,1){50}}
\put(250,40){\line(-1,0){15}}
\put(350,40){\line(1,0){15}}
\put(300,40){\line(0,-1){15}}

\thicklines
\put(225,40){\circle{20}}
\put(375,40){\circle{20}}
\put(300,15){\circle{20}}

\put(250,40){\circle*{5}}
\put(275,40){\circle*{5}}
\put(300,40){\circle*{5}}
\put(325,40){\circle*{5}}
\put(350,40){\circle*{5}}
\put(300,65){\circle*{5}}
\put(300,90){\circle*{5}}

\put(218,37){$-1$}
\put(368,37){$-1$}
\put(297,12){$1$}

\put(320,0){$\widetilde{E}_6$--(2)}

\put(-32,25){1}
\put(-7,25){2}
\put(18,25){3}
\put(43,25){2}
\put(68,25){1}
\put(11,62){2}
\put(11,87){1}
\put(248,25){1}
\put(273,25){2}
\put(291,43){3}
\put(323,25){2}
\put(348,25){1}
\put(291,62){2}
\put(291,87){1}

\end{picture}

\vspace{1cm}

\begin{picture}(400,150)(-50,0)

\thinlines
\put(-30,40){\line(1,0){100}}
\put(20,40){\line(0,1){50}}
\put(-30,40){\line(-1,0){15}}
\put(70,40){\line(1,0){15}}
\put(45,40){\line(0,1){15}}

\thicklines
\put(-55,40){\circle{20}}
\put(95,40){\circle{20}}
\put(45,65){\circle{20}}

\put(-30,40){\circle*{5}}
\put(-5,40){\circle*{5}}
\put(20,40){\circle*{5}}
\put(45,40){\circle*{5}}
\put(70,40){\circle*{5}}
\put(20,65){\circle*{5}}
\put(20,90){\circle*{5}}

\put(-62,37){$-1$}
\put(88,37){$-1$}
\put(42,62){$0$}

\put(10,0){$\widetilde{E}_6$--(3)}

\thinlines
\put(250,40){\line(1,0){100}}
\put(300,40){\line(0,1){50}}
\put(250,40){\line(-1,0){15}}
\put(350,40){\line(1,0){15}}
\put(350,40){\line(0,1){15}}

\thicklines
\put(225,40){\circle{20}}
\put(375,40){\circle{20}}
\put(350,65){\circle{20}}

\put(250,40){\circle*{5}}
\put(275,40){\circle*{5}}
\put(300,40){\circle*{5}}
\put(325,40){\circle*{5}}
\put(350,40){\circle*{5}}
\put(300,65){\circle*{5}}
\put(300,90){\circle*{5}}

\put(218,37){$-1$}
\put(368,37){$-1$}
\put(343,62){$-1$}

\put(290,0){$\widetilde{E}_6$--(4)}

\put(-32,25){1}
\put(-7,25){2}
\put(18,25){3}
\put(43,25){2}
\put(68,25){1}
\put(11,62){2}
\put(11,87){1}
\put(248,25){1}
\put(273,25){2}
\put(298,25){3}
\put(323,25){2}
\put(348,25){1}
\put(291,62){2}
\put(291,87){1}

\end{picture}

\vspace{1cm}

\begin{picture}(400,150)(-50,0)
\thinlines
\put(-30,40){\line(1,0){100}}
\put(20,40){\line(0,1){50}}
\put(-30,40){\line(-1,0){15}}
\put(70,40){\line(1,0){15}}
\put(105,40){\line(1,0){10}}

\thicklines
\put(-55,40){\circle{20}}
\put(95,40){\circle{20}}
\put(125,40){\circle{20}}

\put(-30,40){\circle*{5}}
\put(-5,40){\circle*{5}}
\put(20,40){\circle*{5}}
\put(45,40){\circle*{5}}
\put(70,40){\circle*{5}}
\put(20,65){\circle*{5}}
\put(20,90){\circle*{5}}

\put(-62,37){$-1$}
\put(88,37){$-1$}
\put(118,37){$-2$}

\put(10,0){$\widetilde{E}_6$--(5)}

\thinlines
\put(250,40){\line(1,0){100}}
\put(300,40){\line(0,1){50}}
\put(250,40){\line(-1,0){15}}
\put(325,40){\line(0,1){15}}
\put(325,75){\line(0,1){10}}

\thicklines
\put(225,40){\circle{20}}
\put(325,65){\circle{20}}
\put(325,95){\circle{20}}

\put(250,40){\circle*{5}}
\put(275,40){\circle*{5}}
\put(300,40){\circle*{5}}
\put(325,40){\circle*{5}}
\put(350,40){\circle*{5}}
\put(300,65){\circle*{5}}
\put(300,90){\circle*{5}}

\put(218,37){$-1$}
\put(322,62){$0$}
\put(318,92){$-2$}

\put(290,0){$\widetilde{E}_6$--(6)}

\put(-32,25){1}
\put(-7,25){2}
\put(18,25){3}
\put(43,25){2}
\put(68,25){1}
\put(11,62){2}
\put(11,87){1}
\put(248,25){1}
\put(273,25){2}
\put(298,25){3}
\put(323,25){2}
\put(348,25){1}
\put(291,62){2}
\put(291,87){1}

\end{picture}

\vspace{1.5cm}

\begin{picture}(400,100)(-50,0)

\put(-50,90){\fbox{ Type $\widetilde{E}_7$ }}

\thinlines
\put(50,40){\line(1,0){165}}
\put(125,40){\line(0,1){25}}
\put(50,40){\line(-1,0){15}}
\thicklines
\put(225,40){\circle{20}}
\put(25,40){\circle{20}}

\put(50,40){\circle*{5}}
\put(75,40){\circle*{5}}
\put(100,40){\circle*{5}}
\put(125,40){\circle*{5}}
\put(150,40){\circle*{5}}
\put(175,40){\circle*{5}}
\put(200,40){\circle*{5}}
\put(125,65){\circle*{5}}

\put(218,37){$-1$}
\put(18,37){$-1$}

\put(48,25){1}
\put(73,25){2}
\put(98,25){3}
\put(123,25){4}
\put(148,25){3}
\put(173,25){2}
\put(198,25){1}
\put(130,62){2}

\put(115,0){$\widetilde{E}_7$--(1)}

\end{picture}

\vspace{1cm}

\begin{picture}(400,100)(-50,0)
\thinlines
\put(50,40){\line(1,0){165}}
\put(125,40){\line(0,1){25}}
\put(100,40){\line(0,1){15}}
\thicklines
\put(225,40){\circle{20}}
\put(100,65){\circle{20}}

\put(50,40){\circle*{5}}
\put(75,40){\circle*{5}}
\put(100,40){\circle*{5}}
\put(125,40){\circle*{5}}
\put(150,40){\circle*{5}}
\put(175,40){\circle*{5}}
\put(200,40){\circle*{5}}
\put(125,65){\circle*{5}}

\put(218,37){$-1$}
\put(97,62){$1$}

\put(48,25){1}
\put(73,25){2}
\put(98,25){3}
\put(123,25){4}
\put(148,25){3}
\put(173,25){2}
\put(198,25){1}
\put(130,62){2}

\put(115,0){$\widetilde{E}_7$--(2)}

\end{picture}

\vspace{1cm}

\begin{picture}(400,100)(-50,0)
\thinlines
\put(50,40){\line(1,0){165}}
\put(125,40){\line(0,1){25}}
\put(125,65){\line(0,1){15}}
\thicklines
\put(225,40){\circle{20}}
\put(125,90){\circle{20}}

\put(50,40){\circle*{5}}
\put(75,40){\circle*{5}}
\put(100,40){\circle*{5}}
\put(125,40){\circle*{5}}
\put(150,40){\circle*{5}}
\put(175,40){\circle*{5}}
\put(200,40){\circle*{5}}
\put(125,65){\circle*{5}}

\put(218,37){$-1$}
\put(122,87){$0$}

\put(115,0){$\widetilde{E}_7$--(3)}

\put(48,25){1}
\put(73,25){2}
\put(98,25){3}
\put(123,25){4}
\put(148,25){3}
\put(173,25){2}
\put(198,25){1}
\put(130,62){2}
\end{picture}

\vspace{2cm}

\begin{picture}(400,100)(-50,0)

\put(-50,90){\fbox{ Type $\widetilde{E}_8$ }}

\thinlines
\put(50,40){\line(1,0){175}}
\put(175,40){\line(0,1){25}}
\put(50,40){\line(-1,0){15}}
\thicklines
\put(25,40){\circle{20}}

\put(50,40){\circle*{5}}
\put(75,40){\circle*{5}}
\put(100,40){\circle*{5}}
\put(125,40){\circle*{5}}
\put(150,40){\circle*{5}}
\put(175,40){\circle*{5}}
\put(200,40){\circle*{5}}
\put(225,40){\circle*{5}}
\put(175,65){\circle*{5}}

\put(18,37){$-1$}

\put(48,25){1}
\put(73,25){2}
\put(98,25){3}
\put(123,25){4}
\put(148,25){5}
\put(173,25){6}
\put(198,25){4}
\put(223,25){2}
\put(180,62){3}

\end{picture}
\vspace{0.8cm}


\end{document}